\documentclass[12pt]{article}
\usepackage{amssymb}
\usepackage{amsthm}

\usepackage{extarrows}
\usepackage{cases}
\usepackage{graphicx}
\usepackage{mathrsfs}
\usepackage{caption}
\usepackage{graphics}
\usepackage[misc]{ifsym}
\usepackage{color}
\usepackage{indentfirst}
\usepackage{times}
\usepackage{cite}

\numberwithin{figure}{section}
\numberwithin{equation}{section}

\newtheorem{theorem}{Theorem}[section]
\newtheorem{proposition}[theorem]{Proposition}
\newtheorem{definition}[theorem]{Definition}
\newtheorem{corollary}[theorem]{Corollary}
\newtheorem{lemma}[theorem]{Lemma}
\newtheorem{remark}[theorem]{Remark}
\newtheorem{example}[theorem]{Example}

\DeclareMathOperator*{\limla}{\lim_{\longleftarrow}}

\newcommand{\cA}{{\mathcal A}}

\newcommand{\cC}{{\mathcal C}}
\newcommand{\cD}{{\mathcal D}}

\newcommand{\cH}{{\mathcal H}}

\newcommand{\cK}{{\mathcal K}}
\newcommand{\cL}{{\mathcal L}}

\newcommand{\cN}{{\mathcal N}}

\newcommand{\cW}{{\mathcal W}}

\newcommand{\mB}{{\mbox B}}

\newcommand{\sK}{{\mathscr K}}
\newcommand{\sM}{{\mathscr M}}

\newcommand{\vH}{\check{H}}

\def\R{\mathbb{R}}

\def\N{\mathbb{N}}

\def\1{\mathbb{1}}

\def\disp{\displaystyle}

\def\bc{\begin{center}}
\def\ec{\end{center}}
\def\be{\begin{equation}}
\def\ee{\end{equation}}
\def\ba{\begin{array}}
\def\ea{\end{array}}
\def\benu{\begin{enumerate}}
\def\eenu{\end{enumerate}}
\def\bt{\begin{theorem}}
\def\et{\end{theorem}}
\def\bl{\begin{lemma}}
\def\el{\end{lemma}}
\def\bco{\begin{corollary}}
\def\eco{\end{corollary}}
\def\bn{\begin{numcases}}
\def\en{\end{numcases}}
\def\br{\begin{remark}}
\def\er{\end{remark}}
\def\bd{\begin{definition}}
\def\ed{\end{definition}}
\def\bp{\begin{proposition}}
\def\ep{\end{proposition}}
\def\bo{\begin{proof}}
\def\eo{\end{proof}}
\def\bx{\begin{example}\rm}
\def\ex{\end{example}}

\def\pa{\partial}
\def\al{\alpha}\def\b{\beta}
\def\De{\Delta} \def\de{\delta}
\def\na{\nabla}
\def\lam{\lambda} \def\Lam{\Lambda}

\def\ve{\varepsilon}
\def\sig{\sigma}\def\Sig{\Sigma}\def\vsig{\varsigma}
\def\vp{\varphi}
\def\w{\omega}\def\W{\Omega}
\def\gam{\gamma}\def\Gam{\Gamma}

\def\~{\widetilde}

\def\ol{\overline}

\def\Cap{\bigcap}
\def\Cup{\bigcup}
\def\ra{\rightarrow}

\def\Ra{\Rightarrow}

\def\8{\infty}
\def\X{\times}

\def\mb{\mbox}
\def\di{{\rm d}}

\def\es{\emptyset}
\def\lan{\langle}
\def\ran{\rangle}
\def\sm{\setminus}

\def\ss{\subset}
\def\ssnq{\subsetneq}

\def\Hs{\hspace{0.8cm}}
\def\hs{\hspace{0.4cm}}
\def\Vs{\vskip10pt}
\def\vs{\vskip5pt}

\def\({\left(}
\def\){\right)}

\title{Compactly Generated Shape Index Theory\\ and its Application to a Retarded Nonautonomous Parabolic Equation}

\author{Jintao Wang$^*$\\
{\small Center for Mathematical Sciences \& School of Mathematics and Statistics,}\\
{\small Huazhong University of Science and Technology, Wuhan 430074, China}\\
Jinqiao Duan\\
{\small Department of Applied Mathematics, Illinois Institute of Technology,}\\
{\small Chicago IL 60616, USA}\\
Desheng Li\\
{\small School of Mathematics, Tianjin University, Tianjin 300350, China}}

\date{ }
\begin{document}
\thispagestyle{empty}

\maketitle

\begin{center}
\begin{minipage}{14cm}
\noindent{\bf Abstract}\, We establish the compactly generated shape (H-shape) index theory for local semiflows on complete metric spaces via more general shape index pairs, and define the H-shape cohomology index to develop the Morse equations.
The main advantages are that the quotient space $N/E$ is not necessarily metrizable for the shape index pair $(N,E)$ and $N\sm E$ need not to be a neighborhood of the compact invariant set.
Moreover, in this new theory, the phase space is not required to be separable.
We apply H-shape index theory to an abstract retarded nonautonomous parabolic equation to obtain the existence of bounded full solutions.
\Vs
\Vs
\end{minipage}
\end{center}
\footnote[0]{\vskip-5.5mm
\noindent$^{*}$Corresponding author.\\
E-mail address: wangjt@hust.edu.cn (J.T. Wang); duan@iit.edu (J.Q. Duan); lidsmath@tju.edu.cn (D.S. Li).\\
This work was supported by NSFC grants 11771449 and 11801190.}
\newpage
\setcounter{page}{1}
\section{Introduction}

Index theory plays a significant role in the development of dynamical systems, including topological degree, Morse index, Conley index, Maslov index and shape index.
Among these indices, Conley and shape indices were used to describe the topological property of invariant sets.
Conley index was introduced by Conley \cite{Con1} and extended to local semiflows on complete metric spaces by Rybakowski et al \cite{Ryb} later.
We briefly present the basic idea of Conley index theory below.

Let $K$ be a compact isolated  invariant set of a semiflow $\Phi$ on a complete metric space $X$.
Some appropriate homotopies induced by the semiflow $\Phi$ can help to show that all the pointed quotient spaces $(N/E,[E])$ of Conley index pairs $(N,E)$ have the same homotopy type.
Recall that a Conley index pair is a pair of suitable closed sets $(N,E)$, where $N$ is an isolating neighborhood of $K$, and $E$ is an exit set of $N$.
Then the homotopy Conley index $h(\Phi,K)$ of $K$ is defined to be the homotopy type of the pointed space $(N/E,[E])$.

Shape, invented by Borsuk \cite{Bor} for metric spaces, is a more general topological concept than homotopy type, to describe the topology of spaces with more complicated structures.
Since spaces with the same homotopy type have also the same shape, for the compact isolated invariant set $K$ given above, one can immediately define the shape index $s(K)$ of $K$ as
\be\label{1.1}s(K)={\rm Sh}(N/E,[E]),\ee
where $(N,E)$ is a suitable index pair and ${\rm Sh}$ denotes the shape functor.
This setting is the basic idea of defining shape index of the compact isolated invariant set $K$ in \cite{Wang} and this present paper.

Shape index was first introduced by Robbin and Salamon \cite{Rob} for the flows on a compact smooth manifold.
Their approach to the shape index theory was further developed in the works of Mrozek \cite{Mroz} and S\'{a}nchez-Gabites \cite{Sanc2} for dynamical systems on locally compact spaces.
The case where the phase space is non-locally compact is much more complicated.
Kapitanski and Rodnianski in \cite{Kap} proved that the global attractor of a semiflow on complete metric spaces has a shape of the phase space.
Their work was extended to isolated invariant sets of flows on locally compact metric spaces in \cite{Sanj1} by Sanjurjo, in which the author also considered semiflows on non-locally compact spaces (see \cite{Sanj1}, Section 6).
It was shown that if a semiflow $\Phi$ is two-sided when restricted on the unstable manifold $W^u(K)$, then the shape index of $K$ can be successfully  calculated via its unstable manifold.

For more general cases, the authors in \cite{Wang} used quotient flows to establish the shape Conley index theory via index pairs.
The semiflow $\Phi$ is assumed to be local, asymptotically compact on complete metric spaces and more generally, is not supposed to be two-sided on the unstable manifolds.
The index pair $(N,E)$ of isolated invariant sets used therein is a {\em shape index pair}, different from Conley index pairs.
The shape index pair $(N,E)$ of a compact isolated invariant set $K$ is a pair of closed sets $N$ and $E$ possessing the following properties:\\
\indent(i) $\ol{N\sm E}$ is strongly admissible with $E$ being an exit set of $N$;\\
\indent(ii) $K$ is the maximal compact invariant set in $\ol{N\sm E}$; and\\
\indent(iii) $N\sm E$ contains a local unstable manifold of $K$.\\
To use the Borsuk's shape to give the definition (\ref{1.1}), we need some additional assumptions to guarantee the metrizability of $N/E$ under the quotient topology for a shape index pair $(N,E)$.
It is clear that the shape index pair can be constructed by using local unstable manifolds $W^u_N(K)$ and their appropriate {\em sections}, since the compactness of $W^u_N(K)$ (\cite{Ryb}) can make the quotient space metrizable.
\vs

In this paper, we employ the compactly generated shape defined for general Hausdorff spaces in Rubin and Sanders \cite{Rub}, which allows us to remove the additional assumptions in \cite{Wang}, such as the separability of the phase space and the compactness of the exit set (this was not mentioned clearly in \cite{Wang}).
Then we develop a new type of shape index theory for local semiflows on complete metric spaces via much more general index pairs $(N,E)$.

Since $N/E$ is a normal Hausdorff space for a closed pair $(N,E)$ in a metric space, we adopt the compactly generated shape (H-shape for short, denoted by ${\rm Sh}_{\rm H}$).
We still use the shape index pair $(N,E)$ stated above for a compact isolated invariant set, but with no assumption that the quotient space $N/E$ is metrizable.
Thanks to the consequence in \cite{Gir1} that the compact global attractor and the whole phase space have the same H-shape for semiflows on Hausdorff spaces,
we prove that the pointed spaces $(N/E,[E])$ have the same H-shape for all shape index pairs $(N,E)$, by the similar strategy in \cite{Wang}.
Thus we can define the {\em compactly generated shape index} (H-shape index) $s(K)$ as the H-shape ${\rm Sh}_{\rm H}(N/E,[E])$ of $(N/E,[E])$.
And H-shape index also has the continuation property.

Since the basic idea of establishing shape index theory in this paper resembles that of \cite{Wang}, we necessarily present a precise comparison between them.
Firstly, the paper \cite{Wang} provided the main method and process of quotient flows, by aid of which we developed the shape index theory; in this paper we use the same routine for H-shape index theory.
Secondly, in \cite{Wang}, to define shape Conley index, we are only allowed to choose shape index pairs $(N,E)$ such that $N/E$ is metrizable under the quotient topology,
but in this paper, for H-shape index, we can pick an arbitrary shape index pair without any extra requirement.
Thirdly, due to the lack of metric on the quotient space, the quotient flow used in this paper is merely defined on normal Hausdorff spaces.
Therefore, the relevant results for quotient flows are of new versions and all based on the dynamical systems on topological spaces (see \cite{Li3}).
Moreover, in order to give the Morse equations with respect to H-shape index, we develop somehow cohomological theory and the exactness property for compactly generated shape in this paper.

Note that, H-shape is defined via the direct systems of compact subspaces of the given Hausdorff space $X$ and (ANR-)shape maps between them.
When considering the \v Cech cohomology, which is ANR-shape (see Section 4 below) invariant (\cite{Mar2}), we obtain an inverse system $G^*$ of \v Cech cohomology groups of the compact subspaces of $X$.
We know from Appendix 3.F in Chapter 3 of \cite{Hat} that, the inverse limit of the inverse system $G^*$ may not be isomorphic to the \v Cech cohomology group of $X$,
although the inverse limits are equivalent for equivalent inverse systems.
Thus it remains a question whether the \v Cech cohomology group is H-shape invariant.
However, we can avoid answering this question by using the H-shape cohomology groups, denoted by $C\vH^*$.

Due to the definition of H-shape index of a compact isolated invariant set $K$,
we are allowed to define the cohomology index of $K$ as the H-shape cohomology groups $C\vH^*(N/E,[E])$ for a shape index pair $(N,E)$.
This is sufficient for us to develop the Morse theory.

In our situation, similar to the shape Conley index given in \cite{Wang},
the H-shape index $s(K)$ and Morse equations of an isolated invariant set $K$ can be calculated by using either the Conley index pairs or local unstable manifolds,
which greatly increases the flexibility of calculations.
Moreover, the phase space is not required to be separable.
As an application to illustrate these advantages, we consider an abstract retarded nonautonomous parabolic equation and use the H-shape index to obtain the existence of bounded full solutions.
\vs

This paper is organized as follows.
In Section 2 we present some necessary notions and results in the theory of homotopy on quotient spaces and dynamical system on Hausdorff spaces.
Some necessary results of quotient flows defined on quotient spaces of Wa\.{z}ewski pairs are given in Section 3.
Section 4 is the central part of this paper, in which we introduce the concept of shape index pairs, define compactly generated shape indices for isolated invariant sets and illustrate the continuation property by a simple example.
Section 5 consists of the definition of H-shape cohomology index and the establishment of Morse equations.
In Section 6, we consider an application of H-shape index to an abstract retarded nonautonomous parabolic equation.

\section{Preliminaries}
In this section we collect some necessary notions and results in the theory of topology and dynamical systems on Hausdorff spaces (see \cite{Li3}).
The reader is supposed to be familiar with basic knowledge of algebraic topology.

\subsection{HEP and homotopy equivalence}

Let $X$ be a topological space.
Given a closed subset $A$ of $X$, the  pair $(X,A)$ is said to have the {\em homotopy extension property} (HEP for short), if for every space $Y$ and continuous  mapping $F:\,X\X\{0\}\cup A\X I\ra Y$, there exists a continuous map $\~{F}:\,X\X I\ra Y$ such that $\~{F}$ is an extension of $F$.
\bp[\rm\cite{Ryb}]\label{p2.1}
The pair $(X,A)$ has the HEP if and only if $A$ is a strong deformation retract of one of its open neighborhoods.
\ep
Let $A$ and $B$ be two closed subsets of $X$.
The {\em quotient space} $B/A$ is defined as follows.
If $A\neq\es$, then the space $B/A$ is obtained by collapsing $A$ to a single point $[A]$ in $B\cup A$.
If $A=\es$, we choose a single isolated point $*\notin B$ and define $B/A$ to be the space $B\cup\{*\}$ equipped with the sum topology.
In the latter case we still use the notation  $[A]$ to denote the base point $*$.

We have a homotopy equivalence of quotient spaces as follows, see \cite{Wang}.

\bp\label{p2.2}
Let $A$ and $B$ be two closed subsets of $X$.
Suppose $(X,A)$ has the HEP and that $B$ is a strong deformation retract of $A$.
Then $(X/A,[A])\simeq (X/B,[B])$.
\ep
\subsection{Local semiflows}
Let $X$ be a Hausdorff topological space.
A {\em local semiflow} $\Phi$ on $X$ is a continuous map $\Phi:\,\cD(\Phi)\ra X$, where $\cD(\Phi)$  is an open subset of $\R^+\X X$, and $\Phi$ enjoys the following properties:
\vs (1) for each $x\in X$, there exists $0<T_x\leq\8$ such that $(t,\,x)\in\cD(\Phi)$ if and only if $0\leq t<T_x$;
\vs (2) $\Phi(0,\,x)=x$ for all $x\in X$;
\vs (3) if $(t+s,\,x)\in\cD(\Phi)$, where $t$, $s\in \R^+$, then $\Phi(t+s,\,x)=\Phi(t,\,\Phi(s,\,x))$,\\
where $\cD(\Phi)$ is the domain of $\Phi$.
In the case when $\cD(\Phi)=\R^+\X X$, we simply call $\Phi$  a {\em global semiflow}.

Let $\Phi$ be a given local semiflow on $X$.
For notational convenience, we will rewrite $\Phi(t,\,x)$ as $\Phi(t)x$.

A subset $N$ of $X$ is said to be {\em admissible}, if for arbitrary sequences $x_n\in N$ and $t_n\ra+\8$ with $\Phi([0,t_n])x_n\subset N$ for all $n$, the sequence of the end points $\Phi(t_n)x_n$ has a convergent subsequence.
If additionally $\Phi$ {\it does not explode} in $N$, i.e., we have $T_x=\8$ whenever $\Phi([0,\,T_x))x\ss N$ for all $x\in N$, we say $N$ is {\em strongly admissible}.

Since $X$ may be an infinite-dimensional space, to overcome the difficulty due to the lack of compactness of $X$, we always assume that $\Phi$ is {\em asymptotically compact}, that is, each bounded subset of $X$ is admissible.
It is well known that this condition is naturally satisfied by many important examples from applications, see \cite{Ryb,Tem}.

A {\it solution} (or {\it trajectory}) on an interval $J\subset \R$ is a map $\gam:\,J\ra X$ satisfying
$$
\gam(t)=\Phi(t-s)\gam(s),\Hs\mb{for all }\,s,\,t\in J,\,\,s\leq t.
$$
A {\em full solution} $\gam$ is a solution defined on the whole real line $\R$.
If $x\in X$ is such that $\Phi(t)x=x$ for all $t\geq0$, we say $x$ is an {\em equilibrium}.

The {\em $\w$-limit set} and {\em $\al$-limit set} of a solution $\gam$ are defined as follows.
If $\gam$ is defined on an interval containing $[0,\8)$, it is defined that
$$
\w(\gam)=\{y\in X:\,\mb{there exists }t_n\ra\8\mb{ such that }\gam(t_n)\ra y\}.
$$
If $\gam$ is defined on an interval containing $(-\8,0]$, it is defined that
$$
\al(\gam)=\{y\in X:\,\mb{there exists }t_n\ra-\8\mb{ such that }\gam(t_n)\ra y\}.
$$
For an $x\in X$ with $T_x=\8$, we define $\w(x)=\w(\gam)$ with $\gam(t)=\Phi(t)x$ for $t\geq0$.

A set $A$ is said to be {\em invariant} if $\Phi(t)A=A$ for all $t\geq0$.
For $A\ss X$, we denote by $I(A)$ the maximal invariant set in $A$.
When a closed set $A$ is strongly admissible, one can easily verify $I(A)$ is compact (see Theorem 4.5, Chap. 1 in \cite{Ryb}).

An invariant set $A\ss X$ is said to be {\em isolated}, if $A$ has a neighborhood $N$ such that $A=I(N)$.
Accordingly, a neighborhood $N$ of $A$ such that $A=I(N)$ is called an {\em isolating neighborhood of $A$}.

Given an invariant set $A$ with $A\ss N\ss X$, we define the {\em local stable and unstable manifold}, $W^s_N(A)$ and $W^u_N(A)$ of $A$ in $N$ as follows:
$$
W^s_N(A):\,=\Cup_{\w(\gam)\ss A}\{\gam(t):\,\gam([0,\,\8))\ss N,\,t\in[0,\,\8)\},
$$$$
\mb{and}\hs W^u_N(A):\,=\Cup_{\al(\gam)\ss A}\{\gam(t):\,\gam((-\8,\,0])\ss N,\,t\in(-\8,\,0]\},
$$
where $\gam$ is a solution and $\w(\cdot)$ and $\al(\cdot)$ are limit sets.
If $N=X$ is the whole phase space, we simply write $W^s(A)=W^s_X(A)$ and $W^u(A)=W^u_X(A)$.

\subsection{Attractors}

Here we use the attractor theory of topological spaces stated in \cite{Li3}, which is a generalisation of the attractor theory in metric spaces \cite{Marz,Tem}.

Let $X$ be a Hausdorff topological space and $A,\,B$ be two subsets of $X$.
We say that {\em $A$ attracts $B$}, if $T_x=\8$ for all $x\in B$ and moreover, for an arbitrary neighborhood $U$ of $A$ there exists $T>0$ such that
$$
\Phi(t)B\ss U,\Hs\mb{for all }t>T.
$$
A nonempty sequentially compact invariant set $\cA\ss X$ is said to be an {\em attractor} of $\Phi$, if it attracts a neighborhood $U$ of $\ol{\cA}$ and $\cA$ is the maximal sequentially compact invariant set in $U$.

\br
This definition of attractor differs from the setting in \cite{Li3}, where the authors consider $U$ to be a neighborhood of $A$ if $\ol A\ss{\rm int}U$.
Here we will adopt the concept in common sense that $U$ is a neighborhood of $A$ provided $A\ss{\rm int}U$.
In this sense, this definition of attractor is the same as that in \cite{Li3} in essence.

Here we use sequential compactness over compactness, due to the fact that these two concepts are not equivalent in general topological spaces.
Moreover, we can make good use of the convergence of sequences under sequential compactness (in comparison to \cite{Marz}).

Particularly, for metric spaces, sequential compactness is equivalent to the compactness.
Consequently, if $X$ is metrizable, this definition of attractors is equivalent to those given in \cite{Tem}.
Precisely, $\cA\ss X$ is an attractor of $\Phi$ in $X$, if and only if $\cA$ is nonempty, compact and invariant and attracts a neighborhood of itself.
\er

Let $\cA$ be an attractor. The set $\W(\cA)=\{x\in X:\cA\mb{ attracts }x\}$ is called the {\em region of attraction} (or {\em attraction basin}) of  $\cA$.
One can easily verify that $\W(\cA)$ is open; moreover, $\cA$ attracts each compact subset of $\W(\cA)$, see \cite{Li3}.
In the case when  $\W(\cA)=X$,  we simply call $\cA$  the {\em global attractor} of $\Phi$.

Let $K\ss X$ be a closed subset and $U$ be a subset of $X$ with $K\ss U$.
A continuous function $\zeta:U\ra\R^+$ is called a {\em $\cK_0$ function} of $K$ on $U$, if
$$
\zeta(x)=0\Longleftrightarrow x\in K.
$$
If moreover the {\em level set}
$$
\zeta^a=\{x\in U:\,\zeta(x)\leq a\}
$$
is closed in $X$ for every $a\geq0$, we say $\zeta$ is a {\em $\cK_0^\8$ function} of $K$ on $U$.

If $X$ is a metric space and $A$ is a nonempty closed subset of $X$, then the distance $\di(x,A)$ is a $\cK_0^\8$ function of $A$ on $X$.
If $B$ is another nonempty closed subset of $X$ with $A\cap B=\es$, then the function defined as
$$\frac{\di(x,A)}{\di(x,B)},\hs x\in X\sm B$$
is a $\cK_0^\8$ function of $A$ on $X\sm B$.
Thus we conclude a simple lemma.

\bl\label{l2.6} Let $A$ be a closed subset and $U$ be an open subset of a metric space $X$ with $A\ss U$.
Then, there is a $\cK_0^\8$ function of $A$ on $U$.\el

Let $\cA$ be a closed attractor and $\W:\,=\W(\cA)$ be the region of attraction of $\cA$.
A nonnegative continuous function $\zeta:\,\W\ra\R^+$ is said to be a {\em Lyapunov function} of $\cA$, if $\zeta$ is a $\cK_0$ function of $\cA$ on $\W$, and for $x\in\W\sm\cA$ and $t>0$, we have $\zeta(\Phi(t)x)<\zeta(x)$.
The existence of Lyapunov function for an attractor in Hausdorff spaces is given in the following proposition, a result similar to that for other spaces (\cite{JuQi,Kap,Li1,Li3}).
\bp\label{p:2.14}
Let $X$ be a Hausdorff space.
Assume that the attractor $\cA$ is closed and has a $\cK_0$ function $\psi$ on $\W$.
Then $\cA$ has a Lyapunov function $\zeta$ on $\W$.
What is more, if $\psi$ is $\cK_0^\8$ of $\cA$ on $\W$, so is $\zeta$.
\ep
\bo
The first conclusion comes from Theorem 5.1 in our earlier paper \cite{Li3}.
In order to prove it, we defined for $x\in\W$,
\be\label{3.1}
\psi_1(x)=\sup_{t\geq0}\psi(\Phi(t)x)\hs\mb{and}\hs\zeta(x)=\psi_1(x)+\int_0^{\8}e^{-t}\psi_1(\Phi(t)x)\di t.
\ee
By a standard argument (\cite{Kap,Li3}), we showed that $\psi_1$ is continuous and $\zeta$ is the $\cK_0$ Lyapunov function we want.

Now we assume $\psi$ is $\cK_0^\8$ of $\cA$ on $\W$ and show the second conclusion.
It suffices to show that $\zeta$ defined in (\ref{3.1}) is $\cK_0^\8$ of $\cA$ on $\W$, i.e., for every $a>0$, $\zeta^a$ is closed in $X$.

Indeed, by definition, we have
$$
\zeta(x)\geq\psi_1(x)\geq\psi(x)\hs\mb{ for every }x\in\W.
$$
This means $\zeta^a\ss\psi^a$.
By the continuity of $\zeta$ on $\W$, $\zeta^a$ is a closed subset of $\psi^a$.
The assumption that $\psi$ is $\cK_0^\8$ of $\cA$ on $\W$ implies that $\psi^a$ is closed in $X$.
Then we immediately obtain that $\zeta^a$ is closed in $X$.
The proposition is proved.
\eo

\subsection{Morse Decomposition}

For the reader's convenience, we recall briefly the definition of Morse decompositions of invariant sets for the dynamical systems on topological spaces (see, \cite{Li3} or more classically \cite{Con1,Kap,Li1,Ryb}).

Let $X$ be a topological space and $K$ a compact invariant set.
Then the restriction $\Phi|_K$ of $\Phi$ on $K$ is a semiflow on $K$.
A set $A\ss K$ is called an {\em attractor of $\Phi$ in $K$}, if it is an attractor of $\Phi|_K$.
Note that this attractor in $K$ is a restricted one, defined locally.

Let $A$ be an attractor of $\Phi$ in $K$.
The set $A^*=\{x\in K:\,\w(x)\cap A=\es\}$ is called the {\em repeller} dual to $A$ relative to $K$.
Accordingly,  $(A,\,A^*)$ is called an {\em attractor-repeller pair} in $K$.

Let $K$ be a compact invariant set. An ordered collection
$$
\sM=\{M_1,\,\cdots,\,M_n\}
$$
of subsets $M_k\ss K$ is called a {\em Morse decomposition} of $K$, if there exists an increasing sequence
$\es=A_0\ssnq A_1\ssnq\cdots\ssnq A_n=K$ of attractors in $K$ such that
$$
M_k=A_k\cap A^*_{k-1},\hs 1\leq k\leq n.
$$
The attractor sequence of $A_k$ ($k=0,\,1,\,\cdots,\,n$) is often called the {\em Morse filtration} of $K$, and each $M_k$ is called a
{\em Morse set} of $K$.

\br
It is well known that each Morse set is compact and invariant, and moreover, if $K$ is isolated, then so are the Morse sets $M_k$ (see \cite{Ryb}).
\er

\section{Wa\.{z}ewski Pairs and Quotient Flows}\label{s3}

In this section the phase space $X$ is assumed to be a complete metric space.
Given a subset $N\ss X$, define a function $t_N:\,X\ra\R^+\cup\{\8\}$ as
\be\label{3.0}
t_N(x)=\inf\{t\geq0:\,\,\mb{either $t\geq T_x$\,, or }\Phi(t)x\not\in N\},\Hs\mb{for}\, x\in X.
\ee
Note that  for each $x$, $t_N(x)$ is the maximal time such that $\Phi([0,t_N(x))x\subset N$.

Let $N$ and $E$ be two closed  subsets of $X$.
The subset $E$ is said to be {\em $N$-positively invariant}, if for all $x\in E\cap N$ and $t\geq0$, we have $\Phi([0,t))x\ss E$ whenever $\Phi([0,\,t))x\ss N$.

The subset $E$ is said to be an {\em exit set} of $N$, if
\vs
(1) $E$ is $N$-positively invariant; and
\vs (2) for every $x\in N$ with $t_N(x)<T_x$, there exists $t\leq t_N(x)$ such that $\Phi(t)x\in E$.

\bd A pair of closed subsets $(N,\,E)$ of $X$ is called a {\bf Wa\.{z}ewski pair}, if
\vs
(1) $E$ is an exit set of $N$; and
\vs (2) $\ol{N\setminus E}$ is strongly admissible.
\ed

Let there be given a  Wa\.zewski pair  $(N,\,E)$.
Now we consider the quotient space $N/E$.
For notational simplicity, we denote $[A]=\pi(A)$ for $A\subset N\cup E$, where $\pi:N\cup E\ra N/E$ is the usual quotient map.
Define the {\em quotient flow} $\~\Phi$ of $\Phi$ on $N/E$ as follows:

If $\~{x}=[E]$, then $\~{\Phi}(t)\~x\equiv\~x$ for  $t\in\R^+$; and if $\~{x}=[x]$ for some $x\in\ol{N\setminus E}$, then
$$
\~{\Phi}(t)\~{x}=\left\{
\ba{ll}\,[\Phi(t)x],&\mb{ for }t<t_{\ol{N\setminus E}}(x);\\[1ex]
[E],&\mb{ for }t\geq t_{\ol{N\setminus E}}(x).\ea
\right.
$$
Since $E$ is $N$-positively invariant, it can be easily seen that $\~{\Phi}(t)$ is a well defined semigroup on $N/E$.

Observe that $N/E$ is a normal Hausdorff space.
Lemma 3.2 in \cite{Wang} applies here, and we obtain that $\~\Phi$ is a global semiflow on $N/E$ as follows, see also \cite{Li2}.

\bl\label{l3.1}
The quotient flow $\~{\Phi}$ is continuous on $\R^+\X N/E$ and $N/E$ is strongly admissible.
Moreover, if $I(\ol{N\sm E})\cap E=\es$, the equilibrium $[E]$ is an attractor of $\~{\Phi}$ in $N/E$.
\el

Let $(N,E)$ be a Wa\.zewski pair and $\~\Phi$ be the quotient flow on $N/E$.
Then we have the following conclusions.

\bt\label{th3.7} Every attractor of $\~\Phi$ in $N/E$ is closed.\et
\bo Let $\cA$ be an attractor of $\~\Phi$ in $N/E$ and $\pi$ the quotient map from $N\cup E$ to $N/E$.
If $N\cap E=\es$, the conclusion is obvious.
We only consider the case when $N\cap E\ne\es$.

If $[E]\notin \cA$, we know the restricted map $\pi|_{N\sm E}:N\sm E\ra\pi(N\sm E)$ is homeomorphic.
It follows from the sequential compactness of $\cA$ that $\pi^{-1}(\cA)$ is also sequentially compact in $N\sm E$.
Since $N\sm E$ is a subspace of the metric space $X$, $\pi^{-1}(\cA)$ is compact and hence closed.
Thus $\cA$ is surely closed in $N/E$.

If $[E]\in\cA$, since $\pi$ is a closed map, we only need to show $\pi^{-1}(\cA)$ is closed in $N\cup E$.
It suffices to prove that the limit point $x_0$ of every convergent sequence $x_n\in\pi^{-1}(\cA)$ belongs to $\pi^{-1}(\cA)$.
If $x_0\in E$, we are done.
If $x_0\notin E$, we can assume $x_n\in\pi^{-1}(\cA)\sm E$.
Then via $\pi$, we have $\pi(x_n)\ra\pi(x_0)$ and $\pi(x_n)\in\cA$.
By the sequential compactness of $\cA$, we know $\pi(x_0)\in\cA$ and hence $x_0\in\pi^{-1}(\cA)$.
The proof is complete.\eo

\bt\label{th3.6}
Every attractor $\cA$ of $\~\Phi$ has a $\cK_0^\8$ Lyapunov function on $\W(\cA)$.
\et

\bo
By Theorem \ref{th3.7} and Proposition \ref{p:2.14}, we only need to find a $\cK_0^\8$ function $\psi$ of $\cA$ on $\W(\cA)$.
If $N\cap E=\es$ and $\cA=\{[E]\}$, we have that $\W(\cA)=\{[E]\}$.
The function $\psi([E])=0$ is just what we desire.
Thus we only consider the case when $N\cap E\ne\es$ or $\cA\ne\{[E]\}$.

In order to find the function $\psi$ in this case, we necessarily get back to the Wa\.zewski pair $(N,E)$.
Let $\pi:N\cup E\ra N/E$ be the quotient map.
Denote
$$
U=\pi^{-1}(\W(\cA))\hs\mb{ and }\hs A=\pi^{-1}(\cA).
$$
We see that $U$ is open and $A$ is closed in $N\cup E$.
By Lemma \ref{l2.6}, we obtain a $\cK_0^\8$ function $\de$ of $A$ on $U$.

Now define $\psi:\,\W(\cA)\ra\R^+$ such that
$$
\psi([E])=0\hs\mb{and}\hs\psi(\~x)=\de(x)\mb{ for }\~x=\pi(x)\in\W(\cA)\mb{ with }x\in U\sm E.
$$
The map $\psi$ is well defined.
The continuity of $\psi$ is guaranteed by the properties of $\de$ and $\pi$.
Furthermore, $\psi^a=\pi(\de^a)$ is closed in $N/E$ by the closedness of $\pi$, which indicates that $\psi$ is $\cK_0^\8$ of $\cA$ on $\W(\cA)$.
This completes the proof.
\eo

\section{Compactly Generated Shape Index}
In this section we introduce the shape index pairs and define compactly generated shape indices for isolated invariant sets in metric spaces.

For the presentation hereafter, we first introduce the definition of H-shape of pairs of Hausdorff spaces here, see \cite{Gir1} for the case of Hausdorff spaces and originally \cite{Rub}.

\subsection{H-Shape for Pairs of Hausdorff Spaces}

Shape theory was first introduced by Borsuk \cite{Bor} for metric spaces in 1968, and later Marde\v{s}i\'{c} and Segal \cite{Mar2} gave an extension of Borsuk's shape theory via ANR-systems to include compact Hausdorff topological spaces.
We refer to this definition of shape given by Marde\v si\' c and Segal as {\em ANR-shape}, denoted by ${\rm Sh}_{\rm ANR}$.
In 1974 Rubin and Sanders gave a different extension to the realm of Hausdorff spaces, called ``compactly generated shape'', shortly H-shape, see \cite{Rub}.
The establishment of H-shape theory is based on the ANR-shape theory of compact Hausdorff spaces.

In the following, we introduce the definition of H-shape for pairs of Hausdorff spaces in detail as an extension of H-shape for Hausdorff spaces (\cite{Gir1,Rub}).
\Vs

Let $(X,X_0)$ and $(Y,Y_0)$ be pairs of Hausdorff spaces.
If the pairs satisfy that the relation $Y\ss X$ and $Y_0\ss X_0$, we denote this relation by $(Y,Y_0)\ss(X,X_0)$.
A {\em compact pair} $(X,X_0)$ is a pair with both $X$ and $X_0$ being compact Hausdorff spaces.

Given two compact pairs $(K,K_0)$ and $(L,L_0)$, there is a {\em shape map} $f:(K,K_0)\ra(L,L_0)$, which is defined as follows, see \cite{Gir1,Mar2,Rub}.
The shape map $f$ assigns to every pair $(Q,Q_0)$ having the homotopy type of a CW-complex pair and to every homotopy class $\eta:\,(L,L_0)\ra(Q,Q_0)$, a homotopy class $f(\eta):(K,K_0)\ra(Q,Q_0)$, such that, if $(Q',Q'_0)$ is another pair having the homotopy type of a CW-complex pair and $\eta':(L,L_0)\ra(Q',Q'_0)$ is a homotopy class, then if $\xi:(Q,Q_0)\ra(Q',Q'_0)$ is a homotopy class, the commutativity (up to homotopy) of

\begin{picture}(0,80)
\put(30,60){\vector(1,0){70}}
\put(28,55){\vector(3,-2){70}}
\put(118,52){\vector(0,-1){42}}
\put(58,63){\small$\eta$}
\put(-5,57){\small$(L,L_0)$}
\put(103,57){\small$(Q,Q_0)$}
\put(50,25){\small$\eta'$}
\put(103,0){\small$(Q',Q'_0)$}
\put(120,25){\small$\xi$}
\put(140,25){ implies that of }
\put(225,60){\vector(1,0){70}}
\put(223,55){\vector(3,-2){70}}
\put(313,52){\vector(0,-1){42}}
\put(253,63){\small$f(\eta)$}
\put(187,57){\small$(K,K_0)$}
\put(230,25){\small$f(\eta')$}
\put(298,57){\small$(Q,Q_0)$}
\put(298,0){\small$(Q',Q'_0)$}
\put(315,25){\small$\xi$}
\end{picture}

Let $\Lam$ be a directed set.
A {\em CS$^2$-system} is a direct system $X^*=\{(X_\lam,X_{0\lam}),p_{\lam\lam'},\Lam\}$ in the compact shape category of pairs of Hausdorff spaces (see \cite{Mar2}), that is to say, each $(X_\lam,X_{0\lam})$ is a compact pair and if $\lam\leq\lam'$ in $\Lam$, then $p_{\lam\lam'}:(X_\lam,X_{0\lam})\ra(X_{\lam'},X_{0\lam'})$ is a shape map such that
\benu\item[(i)] $p_{\lam\lam}=1_{(X_\lam,X_{0\lam})}$ is the identity shape map,
\item[(ii)] if $\lam\leq\lam'\leq\lam''$, then $p_{\lam'\lam''}p_{\lam\lam'}=p_{\lam\lam''}$.\eenu
A {\em CS$^2$-morphism} $F:X^*\ra Y^*=\{(Y_\mu,Y_{0\mu}),q_{\mu\mu'},M\}$ is a pair $F=(f_\lam,f)$ consisting of an increasing function $f:\Lam\ra M$ and a collection of shape maps $f_\lam:(X_\lam,X_{0\lam})\ra(Y_{f(\lam)},Y_{0f(\lam)})$ such that if $\lam\leq\lam'$ then $q_{f(\lam)f(\lam')}f_\lam=f_{\lam'}p_{\lam\lam'}$, that is to say, the following diagram commutes.

\begin{picture}(0,80)
\put(70,70){\vector(1,0){70}}
\put(70,13){\vector(1,0){70}}
\put(43,63){\vector(0,-1){43}}
\put(168,63){\vector(0,-1){43}}
\put(108,73){\small$f_\lam$}
\put(108,4){\small$f_{\lam'}$}
\put(19,67){\small$(X_\lam,X_{0\lam})$}
\put(143,67){\small$(Y_{f(\lam)},Y_{0f(\lam)})$}
\put(19,10){\small$(X_{\lam'},X_{0\lam'})$}
\put(143,10){\small$(Y_{f(\lam')},Y_{0f(\lam')})$}
\put(25,38){\small$p_{\lam\lam'}$}
\put(170,38){\small$q_{f(\lam)f(\lam')}$}
\end{picture}

Defining the identity $1_{X^*}$ and compositions in the usual way, we finally have a category of CS$^2$-systems and CS$^2$-morphisms between them, denoted by CS$^2$.

Two CS$^2$-morphisms $F,G:X^*\ra Y^*$ are {\em homotopic}, $F\simeq G$, if for each $\lam\in\Lam$, there is $\mu\in M$ with $f(\lam),\,g(\lam)\leq\mu$ such that $q_{f(\lam)\mu}f_\lam=q_{g(\lam)\mu}g_\lam$, i.e., the following commutative diagram.

\begin{picture}(0,90)
\put(32,52){\vector(3,2){35}}
\put(32,35){\vector(3,-2){35}}
\put(93,75){\vector(3,-2){35}}
\put(93,11){\vector(3,2){35}}
\put(48,79){\small$(Y_f(\lam),Y_{0f(\lam)})$}
\put(48,3){\small$(Y_{g(\lam)},Y_{0g(\lam)})$}
\put(39,65){\small$f_\lam$}
\put(112,65){\small$q_{f(\lam)\mu}$}
\put(39,20){\small$g_\lam$}
\put(112,19){\small$q_{g(\lam)\mu}$}
\put(0,40){\small$(X_\lam,X_{0\lam})$}
\put(124,40){\small$(Y_\mu,Y_{0\mu})$}
\end{picture}

Surely the homotopy relation $\simeq$ is a morphism equivalence, see \cite{Rub}.
We say $X^*$ and $Y^*$ have the same {\em homotopy type}, provided there are CS$^2$-morphisms $F:X^*\ra Y^*$ and $G:Y^*\ra X^*$ such that $GF\simeq 1_{X^*}$ and $FG\simeq 1_{Y^*}$;
and we say $F$ is a {\em homotopy equivalence} from $X^*$ to $Y^*$.
\vs

Given a pair $(X,X_0)$ of Hausdorff spaces, let $c(X,X_0)$ be the set of all compact pairs of $(K,K_0)\ss(X,X_0)$ ordered by inclusions, which makes $c(X,X_0)$ a directed set.
Then one has a CS$^2$-system
$$
C(X,X_0)=\{(K,K_0),i_{(K,K_0)(K',K'_0)},c(X,X_0)\}
$$
such that $(K,K_0)\in c(X,X_0)$ and if $(K,K_0)\ss(K',K'_0)$ then $i_{(K,K_0)(K',K'_0)}$ is the inclusion shape map.

\bd
Let $(X,X_0)$ and $(Y,Y_0)$ be pairs of Hausdorff spaces.
If $C(X,X_0)$ and $C(Y,Y_0)$ have the same homotopy type, we say $(X,X_0)$ and $(Y,Y_0)$ have the same {\bf shape}, denoted by ${\rm Sh}_{\rm H}(X,X_0)={\rm Sh}_{\rm H}(Y,Y_0)$.
\ed
\br
Since ANR-shape is defined by the inverse systems of neighborhoods of a given metric space (compact Hausdorff space) in an ANR and the homotopy classes between them, it mainly describes the space from outside, see \cite{Kap,Mar2}.
Here {\em ANR} means the absolute neighborhood retract of metric spaces (or compact Hausdorff spaces) (see \cite{Mar2}).

By comparison, H-shape defined above is via the direct systems of compact subsets of a given Hausdorff space and the (ANR-)shape maps between them (see also \cite{Gir1,Rub}).
Correspondingly, H-shape provides an inner description of the Hausdorff space.
In spite of the definitions in different means and distinct descriptions of the space, ANR-shape and H-shape coincide for compact Hausdorff spaces (\cite{Rub}).
\er
If $\vp:(X,X_0)\ra(Y,Y_0)$ is a continuous map, let $f:c(X,X_0)\ra c(Y,Y_0)$ and $f_{(K,K_0)}:(K,K_0)\ra(\vp(K),\vp(K_0))$ such that $f((K,K_0))=(\vp(K),\vp(K_0))$, which is increasing, and $f_{(K,K_0)}$ is the shape map induced by $\vp|_{(K,K_0)}:(K,K_0)\ra(\vp(K),\vp(K_0))$.
Thus we obtain a CS$^2$-morphism $F=(f_{(K,K_0)},f)$ induced by $\vp$.
If $F$ is a homotopy equivalence from $C(X,X_0)$ to $C(Y,Y_0)$, we say that $\vp$ {\em induces an H-shape equivalence}.

We can see in a straightforward way that H-shape is a homotopy invariant for pairs of Hausdorff spaces, as stated in Rubin and Sanders \cite{Rub}, i.e.,
$$
(X,X_0)\simeq(Y,Y_0)\hs\Ra\hs{\rm Sh}_{\rm H}(X,X_0)={\rm Sh}_{\rm H}(Y,Y_0).
$$
\vs

The following result is a pointed version of Theorem 4.2 in \cite{Gir1}, and the similar results in metric spaces can be found  in \cite{Gir2,Kap,Sanj2}.
\bt\label{t3.2}
Let $X$ be a Hausdorff space, and $\Phi$ be a global semiflow on $X$.
Suppose that $\Phi$ has a compact global attractor $\cA$, and  that the system has an equilibrium $e\in\cA$.
Then the inclusion $(\cA,\,e)\ra(X,\,e)$ induces an H-shape equivalence.
\et

Let $(X,x_0)$ and $(Y,y_0)$ be two pointed spaces.
The {\em wedge sum} $(X,x_0)\vee(Y,y_0)$ and {\em smash product} $(X,x_0)\wedge(Y,y_0)$ are defined, respectively, as follows,
$$
(X,x_0)\vee(Y,y_0)=(\cW,(x_0,y_0)),\hs (X,x_0)\wedge(Y,y_0)=((X\X Y)/\cW,[\cW]),
$$
where $\cW=X\X\{y_0\}\cup\{x_0\}\X Y$.
Similar to the definition for homotopy type (see \cite{Hat,Ryb}), the operations ``$\vee$'' and ``$\wedge$'' can be also defined to H-shape of pointed spaces as follows.
This is also a natural generalisation of the definitions of wedge sum and smash product for H-shape of Hausdorff spaces (see \cite{Rub}).

For pointed spaces $(X,x_0)$, $(X',x'_0)$, $(Y,y_0)$ and $(Y',y'_0)$, if
$${\rm Sh}_{\rm H}(X,x_0)={\rm Sh}_{\rm H}(X',x'_0)\hs\mb{and}\hs{\rm Sh}_{\rm H}(Y,y_0)={\rm Sh}_{\rm H}(Y',y'_0),$$
then by the definition of H-shape of pointed spaces and using the method in \cite{Rub}, we similarly obtain that
$${\rm Sh}_{\rm H}((X,x_0)\vee(Y,y_0))={\rm Sh}_{\rm H}((X',x'_0)\vee(Y',y'_0))$$
$$\mb{and}\hs{\rm Sh}_{\rm H}((X,x_0)\wedge(Y,y_0))={\rm Sh}_{\rm H}((X',x'_0)\wedge(Y',y'_0)).$$
This allows us to define the operators $\vee$ and $\wedge$ for H-shape of pointed spaces as follows,
$$
{\rm Sh}_{\rm H}(X,x_0)\vee{\rm Sh}_{\rm H}(Y,y_0)={\rm Sh}_{\rm H}((X,x_0)\vee(Y,y_0)),
$$$$
\mb{and}\hs{\rm Sh}_{\rm H}(X,x_0)\wedge{\rm Sh}_{\rm H}(Y,y_0)={\rm Sh}_{\rm H}((X,x_0)\wedge(Y,y_0)).
$$

\subsection{H-Shape Index}
From now on we let $X$ be a complete metric space and $\Phi$ a local semiflow on $X$.

\bd\label{d:3.7}
Let $K\ss X$ be a compact isolated  invariant set of $\Phi$.
A Wa\.zewski pair $(N,\,E)$ is said to be a {\bf shape index pair} of $K$, if
\vs
(1) there is a closed admissible neighborhood $U$ of $K$ such that $W_U^u(K)\ss\ol{N\sm E}$;
\vs (2) $K\cap E=\es$; and
\vs (3) $K=I(\ol{N\sm E})$.
\ed

\br (1) A Conley index pair $(N,E)$ of a compact isolated invariant set $K$ is a Wa\.zewski pair such that $N\sm E$ is an isolating neighborhood of $K$ (see \cite{Con1,Ryb}).
The set $\ol{N\sm E}$ surely contains a local unstable manifold of $K$ and $K\cap E=\es$.
Therefore, a Conley index pair is naturally a shape index pair in Definition \ref{d:3.7}.

(2) The shape index pairs $(N,E)$ given in \cite{Wang} are also specific examples of the shape index pairs defined above, since besides the conditions (1) (2) (3) in Definition \ref{d:3.7}, it is also implicated therein that the quotient space $N/E$ is metrizable in quotient topology.
\er

We are now prepared to define the H-shape index via shape index pairs.

\bd
Let $(N,\,E)$ be a shape index pair of $K$.
The {\bf H-shape index} $s(\Phi,K)$ of $K$ is defined as
$$s(\Phi,\,K)={\rm Sh}_{\rm H}(N/E,\,[E]).$$
\ed
When the semiflow $\Phi$ is clear, we will simply write $s(\Phi,K)$ as $s(K)$.
\bx\label{x4.4} (1) If the compact isolated invariant set $K=\es$, we can take $(\es,\es)$ to be the shape index pair of $\es$, and thus we have $s(\es)=\ol0$, where $\ol0$ is the H-shape of a pointed singleton.
Applying this fact, we can determine that $K\ne\es$, if $s(K)\ne\ol0$.
This property is analogous to Conley index.
\vs

(2) By Theorem \ref{t3.2}, if $X$ is a complete metric space and $\cA$ is the global attractor, it is clear that
$$
s(\cA)={\rm Sh}_{\rm H}(X\cup\{*\},*),
$$
where $*\notin X$ is a single isolated point and $X\cup\{*\}$ is endowed the sum topology.
Particularly, if $X$ is a normed linear space, $s(\cA)=\Sig^0$, since $X$ is contractible.
Here and in the sequel, we use $\Sig^n$ to denote the H-shape of pointed $n$-dimensional sphere.
\ex

The following result indicates that the H-shape index is well defined.
\bt[Main Theorem]\label{t:4.1}
The H-shape index $s(K)$ of $K$ is independent of the choice of shape index pairs.
\et

To show Theorem \ref{t:4.1}, we need the following lemma (we omit the proof), which is a new version of Lemma 4.6 and 4.7 in \cite{Wang}, in the framework of topological dynamical systems (given in Section 2) on the quotient space.

\bl\label{l3.3}
Let $K$ be a compact isolated invariant set with its shape index pair $(N,E)$ and let $\~\Phi$ be the quotient flow on $N/E$.
Then $\~\Phi$ has a compact global attractor $\cA$ such that
$$
\cA=W^u([K])\cup\{[E]\}\hs\mb{and}\hs(\cA,[E])\simeq(W_N^u(K)/E,[E]).
$$
\el

The proof of Theorem \ref{t:4.1} is indeed a modification of that of Theorem 4.5 in \cite{Wang}.
For the reader's convenience, we present the main sketch of the proof as follows.
Nevertheless, we need to keep in mind that the referred conclusions in \cite{Wang} involved in the following proof is under the framework of the shape index pair and the quotient flow defined in this paper.
And the proofs of the referred conclusions in \cite{Wang} also work well here, with the Lyapunov functions used therein following from Theorem \ref{th3.6}.
\Vs

\noindent{\it Proof of Theorem \ref{t:4.1}.}
In the case when $W^u(K)=K$, for every shape index pair $(N,E)$, $W^u_N(K)=K$.
By Lemma \ref{l3.3}, the quotient flow on $N/E$ has a global attractor $\cA$.
Moreover, we have $\cA=[K]\cup\{[E]\}$, and so
$$
{\rm Sh}_{\rm H}(N/E,[E])={\rm Sh}_{\rm H}(\cA,[E])={\rm Sh}_{\rm H}([K]\cup \{[E]\},\,[E])={\rm Sh}_{\rm H}(K\cup\{*\},*),
$$
which implies the conclusion of Theorem \ref{t:4.1}.
\vs

Now we only consider the case when $W^u(K)\ne K$.
In this case, for all shape index pairs $(N,E)$ of $K$, we infer that $W^u_{N}(K)\cap E\ne \es$ (see Lemma 4.8 in \cite{Wang}).

Let $(N_1,E_1)$ and $(N_2,E_2)$ be two shape index pairs of $K$. Let $N=N_1\cap N_2$ and $E=E_1\cup E_2$.
Then we have that $(N,E)$ is also a shape index pair of $K$.
We aim to show that $(N_1/E_1,[E_1])$ and $(N_2/E_2,[E_2])$ have the same H-shape, for which we only need to prove that
\be\label{4.2}
{\rm Sh}_{\rm H}(N_k/E_k,[E_k])={\rm Sh}_{\rm H}(N/E,[E]),\hs k=1,\,2.
\ee
In the following we only need to show that (\ref{4.2}) holds true for $k=1$, since the same argument also works for $k=2$.

Set $E^u =W^u_{N }(K)\cap E ,$ and define
$$
\ba{ll}\Gam_1=\{x\in N_1:&\mb{there exist }t\geq0\mb{ and }y\in E^u \mb{ such that }\\
&\Phi([0,\,t])y\ss N_1,
\mb{ and }\Phi(t)y=x\}.\ea
$$
It is easy to check that $E^u $ and $\Gam_1$ are $N$-positively invariant and $N_1$-positively invariant, respectively.
Moreover, we have
\be\label{4.3}
W^u_N(K)\cup \Gam_1=W^u_{N_1}(K).
\ee
Applying Lemma 4.10 in \cite{Wang}, we have an open neighborhood $U$ of $K$ such that $\Gam_1\cap U=\es$.

By using Lemma 4.9 in \cite{Wang}, there is a closed neighborhood $F$ of $E$ in $N\cup E$ with $K\cap F=\es$ such that $(N,F)$ is a shape index pair and has HEP; moreover,
\be\label{4.4}
W^u_N(K)\sm F\ss U,\hs (N/E,[E])\simeq (N/F,[F]).
\ee
Hence ${\rm Sh}_{\rm H}(N/E,[E])={\rm Sh}_{\rm H}(N/F,[F])$.
On the other hand, by Theorem \ref{t3.2}, we deduce
$$
{\rm Sh}_{\rm H}(N/F,\,[F])={\rm Sh}_{\rm H}(\cA',\,[F]),
$$
where $\cA'$ is the global attractor of the quotient flow on $N/F$.
Therefore by Lemma \ref{l3.3} we find that
\be\label{4.5}
\ba{ll}&{\rm Sh}_{\rm H}(N/E,[E])={\rm Sh}_{\rm H}(\cA',[F])\\[1ex]
=&{\rm Sh}_{\rm H}(W^u_N(K)/F,[F])={\rm Sh}_{\rm H}(W^u_N(K)/F^u,[F^u]),\ea
\ee
where $F^u=W^u_N(K)\cap F$.

Let $\Gam=F^u\cup \Gam_1$.
Based on the fact that $\Gam_1\cap U=\es$ and (\ref{4.4}), we have
$$
\ba{ll}&W_N^u(K)\cap\Gam=F^u\cup(W_N^u(K)\cap\Gam_1)\\
=&F^u\cup[W_N^u(K)\cap(F\cup F^c)\cap\Gam_1]=F^u\cup(W_N^u(K)\cap F\cap\Gam_1)=F^u,\ea
$$
where $F^c=X\sm F$. Then we have
$$\ba{ll}
&\(W^u_N(K)/F^u,\,[F^u]\)\cong \(W^u_N(K)/\Gam,\,[\Gam]\)\\[1ex]
\cong&\((W^u_N(K)\cup\Gam)/\Gam,\,[\Gam]\)\cong (\mb{by }(\ref{4.3}))\cong \(W^u_{N_1}(K)/\Gam,\,[\Gam]\).\ea
$$
Therefore by (\ref{4.5}) we obtain that
\be\label{4.6}
{\rm Sh}_{\rm H}(N/E,\,[E])={\rm Sh}_{\rm H}(W^u_{N_1}(K)/\Gam,\,[\Gam]).
\ee

Let $E^u_1=W^u_{N_1}\cap E_1$.
Consider the quotient space $W^u_{N_1}(K)/E_1^u$ along with the quotient flow $\~\Phi_1$.
Let $\pi$ be the quotient map from $W^u_{N_1}(K)$ to $W^u_{N_1}(K)/E_1^u$.
It is clear that
\be\label{4.7}
\(W^u_{N_1}(K)/\Gam,[\Gam]\)\cong \(\pi(W^u_{N_1}(K))/\pi(\Gam),\,[\pi(\Gam)]\).
\ee
Theorem \ref{l3.1} asserts that $[E_1^u]$ is an attractor of $\~\Phi_1$.
Since $\pi(\Gam)$ is positively invariant and contained in the attraction basin of $[E_1^u]$ with $[E_1^u]\in\pi(\Gam)$, then $[E_1^u]$ is a strong deformation retract of $\pi(\Gam)$ (see Proposition 2.5 in \cite{Wang}).
Because $(N,F)$ has HEP, by Proposition \ref{p2.1}, it is easy to see that $(W_N^u(K),F^u)$ has HEP.
Consequently, $(W_{N_1}^u(K),\,\Gam)$ and $\(\pi(W^u_{N_1}(K)),\pi(\Gam)\)$ have HEP as well.
Therefore by Proposition \ref{p2.2} we have
$$\ba{ll}&\(\pi(W^u_{N_1}(K))/\pi(\Gam),[\pi(\Gam)]\)\simeq\(\pi(W^u_{N_1}(K))/\pi(E^u_1),[\pi(E^u_1)]\)\\[1ex]
\cong&\(W^u_{N_1}(K)/E^u_1,[E^u_1]\)\simeq\(W^u_{N_1}(K)/E_1,[E_1]\).\ea
$$
By (\ref{4.6}) and (\ref{4.7}), we have that
$$
\ba{ll}&{\rm Sh}_{\rm H}(N/E,[E])={\rm Sh}_{\rm H}(W^u_{N_1}(K)/E_1,[E_1])\\[1ex]
=&(\mb{by Lemma \ref{l3.3}})={\rm Sh}_{\rm H}(\cA_1,[E_1]),\ea
$$
where $\cA_1$ is the global attractor of the quotient flow on $N_1/E_1$.
Furthermore, by Theorem \ref{t3.2} we conclude that
$$
{\rm Sh}_{\rm H}(N/E,[E])={\rm Sh}_{\rm H}(\cA_1,\,[E_1])={\rm Sh}_{\rm H}(N_1/E_1,\,[E_1]).
$$
The proof is finished now.
\qed

\subsection{Continuation Property}
\label{s4.3}
Similar to the Conley index (see \cite{Ryb}) and shape Conley index, H-shape index also has the continuation property, which involves a continuous family of local semiflows.
We follow the basic concepts given for Conley index theory in \cite{Ryb} below.

Let $X$ be a complete metric space.
For a sequence of local semiflows $\Phi_n$ on $X$, we write $\Phi_n\ra\Phi_0$, if for all sequences $x_n\in X$ and $t_n\in\R^+$ with $x_n\ra x_0$ and $t_n\ra t_0$,
$\Phi_n(t_n)x_n\ra\Phi_0(t_0)x_0$.

Let $\Phi_n$ be a sequence of local semiflows on $X$.
A set $N\ss X$ is said to be {\em strongly $\{\Phi_n\}$-admissible} if for two arbitrary sequences $x_n\in X$ and $0\leq t_n\ra\8$ satisfying $\Phi_n([0,t_n])x_n\ss N$ for all $n\in\N^+$, the sequence of endpoints $\Phi_n(t_n)x_n$ has a convergent subsequence and furthermore, $\Phi_n$ does not explode in $N$ for every $n\in\N^+$.

Let $\Lam$ be a metric space. We write $\Phi_\lam\ra\Phi_{\lam_0}$, if $\Phi_{\lam_n}\ra\Phi_{\lam_0}$ for every sequence $\lam_n\in\Lam$ with $\lam_n\ra\lam_0$.
A {\em continuous family of local semiflows} $\Phi_\lam$ on $X$ is a family of local semiflows such that $\Phi_{\lam}\ra\Phi_{\lam_0}$ for each $\lam_0\in\Lam$ with $\lam\ra\lam_0$.

Now we recall the definition of S-continuity.
The pair $(\Phi_\lam,\,K_\lam)$ is said to be {\em S-continuous at $\lam_0\in\Lam$},
if there is $\de>0$ and a closed subset $N$ of $X$ such that the following two conditions are fulfilled:
\benu\item[(1)] for every $\Phi_\lam$ and $\lam$ with $\rho(\lam,\lam_0)<\de$, the subset $N$ is a strongly admissible closed neighborhood of $K_\lam$;
\item[(2)] Whenever $\lam_n\ra\lam_0$, then $\Phi_{\lam_n}\ra\Phi_{\lam_0}$ and $N$ is $\{\Phi_{\lam_n}\}$-admissible.\eenu
If $(\Phi_\lam,\,K_\lam)$ is S-continuous at each point $\lam\in\Lam$, $(\Phi_\lam,\,K_\lam)$ is said to be {\em S-continuous on $\Lam$}.

Since the Conley index pairs are also shape index pairs for an isolated compact set, by similar discussion for Conley index (see Theorem 12.2 in \cite{Ryb}), one immediately concludes that H-shape index possesses the following property, which is just what we call the {\em continuation property}.
\bt\label{th4.9}
Let $K_\lam$ be a compact isolated invariant set of $\Phi_{\lam}$ for each $\lam$ lying in a connected component $\Lam_0$ of $\Lam$.
Suppose $(\Phi_\lam,\,K_\lam)$ is S-continuous on $\Lam_0$.
Then $s(\Phi_{\lam},K_{\lam})$ is constant for $\lam\in\Lam_0$.\et

\subsection{An Example}
For a better understanding of continuation property and the shape index pair, we consider the initial-boundary problem of the equation
\be\label{4.8}
\left\{\ba{ll}\disp\frac{\pa u}{\pa t}-\De u=\b u(1-u^{p-2}),\hs &x\in\W,\,t>0,\\[1ex]
u(x,t)=0,&x\in\pa\W,\,t\geq0,\\[1ex]
u(x,0)=u_0(x),&x\in\ol\W,\ea\right.
\ee
where $\W\subset \R^n$ is a bounded domain, $p>2$ and $\b>0$ are constants.
Such a problem has an invariant domain
$$
X=\{u\in L^\8(\W)|\,\,0\leq u\leq 1\}.
$$
Hence we have a semiflow on $X$.
Note that $X$ is not separable.
Now we compute the H-shape index of $0$ in $X$.

Let $V=H_0^1(\W)$.
The inner product $\lan\cdot,\cdot\ran$ on $V$ is defined by
$$
\lan u,v\ran=\int_\W\na u\na v\di x,\hs u,\,v\in V,
$$
and the corresponding norm $\|\cdot\|$.
We know that the equation (\ref{4.8}) has a {\em weak solution}
$$
u\in L^{\8}(0,T;X)\cap L^2(0,T;V),
$$
such that
$$
\int_{\W}\frac{\pa u}{\pa t}v\di x+\lan u,v\ran=\int_{\W}\b u(1-u^{p-2})v\di x,
$$
for every $v\in V$ (see \cite{Tem}).

We first consider the case when the phase space is $V\cap X$.
To compute the shape index of $0$ for (\ref{4.8}), we consider the following equation related to (\ref{4.8}),
\be\label{4.9}
\left\{\ba{ll}\disp\frac{\pa u}{\pa t}-(\De+\b)u=-\lam\b u^{p-1},&x\in\W,\,t>0,\\[1ex]
u(x,t)=0,&x\in\pa\W,\,t\geq0,\ea\right.
\ee
where $\lam\in[0,1]$. When $\lam=0$, (\ref{4.9}) is the linearisation of (\ref{4.8}) at $u=0$; when $\lam=1$, (\ref{4.9}) is just (\ref{4.8}).
It is well known that the linear operator $-(\De+\b)$ has only finitely many negative eigenvalues.
Denote the eigenvalues of $-\De$ in $V$ by
$$
0<\mu_1<\mu_2<\cdots<\mu_k<\cdots\ra\8.
$$

Let $\Phi_\lam$ be the semiflow generated by (\ref{4.9}).
We first consider the case when $\b\ne\mu_k$ for all $k\in\N^+$.
By the standard argument (see \cite{Ryb}), the family of semiflows $\Phi_\lam$ is continuous in $\lam\in[0,1]$, and $(\Phi_\lam,\{0\})$ is S-continuous on $[0,1]$.
By the continuation property (Theorem \ref{th4.9}), we have
\be\label{4.10}
s(\Phi_1,\{0\})=s(\Phi_0,\{0\})=\left\{\ba{ll}\ol0,&\b>\mu_1,\,\b\ne\mu_k,\\\Sig^0,&0<\b<\mu_1.\ea\right.
\ee
The second equality in (\ref{4.10}) results from the following discussion.

Since $\b\ne\mu_k$, we let $V_0$ be the subspace of $V$ spanned by the eigenfunctions of $\mu_k$ with $0<\mu_k<\b$.
Pick $R\in(0,1)$ and let
$$N:=\{u\in V_0\cap L^\8(\W):0\leq u\leq R\}\hs\mb{and}\hs E:=\{u\in N:u=R\}.$$
We know that $N$ is a local unstable manifold of $0$ in $V\cap X$ and $E$ is the exit set of $N$.
Hence the pair $(N,E)$ is a shape index pair of $\{0\}$ for $\Phi_0$.
When $0<\b<\mu_1$, we see that $(N,E)=(\{0\},\es)$ and so $s(\Phi_0,\{0\})=\Sig^0$.
When $\b>\mu_1$, the space $N/E$ is contractible and then $s(\Phi_0,\{0\})=\ol0$.

If $\b=\mu_k$ for some $k\in\N^+$, we need to consider the centre manifold of $0$.
It is well known (see Chapter 6, \cite{Hen}) that $0$ is asymptotically stable in its centre manifold.
As a consequence, we obtain that
\be\label{4.11}
s(\Phi_0,\{0\})=\left\{\ba{ll}\ol0,&\b=\mu_k,\,k>1,\\\Sig^0,&\b=\mu_1.\ea\right.
\ee

Now we consider the case when the phase space is $X$.
Let $\Phi$ be the semiflow generated by (\ref{4.8}).
Since there is a natural embedding from $L^\8(\W)$ into $L^p(\W)$, then $\Phi(t)u_0\in V$ for all $u_0\in X$.
The unstable manifold of $0$ in $V\cap X$ is hence also that of $0$ in $X$.
Note that the local unstable manifold $W^u_U(0)$ of $0$ is finite-dimensional, for some closed neighborhood $U$ of $0$ in $X$.
The topologies of $W^u_U(0)$ induced by $X$ and $V\cap X$ are equivalent.
As a result, by summarising (\ref{4.10}) and (\ref{4.11}), we conclude the following result
$$
s(\Phi,\{0\})=\left\{\ba{ll}\ol0,&\b>\mu_1,\\\Sig^0,&0<\b\leq\mu_1.\ea\right.
$$

\br
If we replace the nonlinear term $\b u(1-u^{p-2})$ in (\ref{4.8}) by $\b u(1-|u|^{p-2})$, we are allowed to set the phase space to be $X=\{u\in L^\8(\W)|-1<u<1\}$.
With almost the same argument, we obtain that $s(\{0\})=\Sig^{r}$, where $r$ is defined as
$$r=r(\b):=\sum_{0<\mu_k<\b}r_k,$$
where $r_k$ denotes the multiplicity of the eigenvalue $\mu_k$.
\er
\section{Establishment of Morse Equations}

In this section we study the Morse equations for a compact isolated invariant set associated to H-shape.
To do this, we need to consider some H-shape invariant for these invariant sets.
In our theory, we need the H-shape invariant to possess the following property:\\
\noindent({\bf P}) The pair $(N,E)$ and its quotient space $(N/E,[E])$ have the same H-shape invariant.

By Theorem 2.3 of \cite{Sand1} and Theorem 2.9 of \cite{Sand2}, we know the shape groups (obtained from direct limits of homotopy groups) are H-shape invariants for the pointed Hausdorff spaces.
However, homotopy groups do not meet (P), let alone the shape groups.
Note that the \v Cech cohomology groups satisfy the property (P) and the exactness property (\cite{Eil}).
But we do not know whether the \v Cech cohomology group is H-shape invariant.
Hence it is necessary to develop some new type of cohomology theory based on \v Cech cohomology groups.

\subsection{H-shape Cohomology Groups and Indices}

Given a pair $(X,X_0)$ of Hausdorff spaces, an abelian group $\mathfrak{A}$ and each $q\in\N$, there exists a \v Cech cohomology group $\check H^q(X,X_0;\mathfrak{A})$, which is also abelian.
It is found from the third section of Chapter II in \cite{Mar2} that \v Cech cohomology groups are shape (by the inverse systems) invariant for pairs of topological spaces, including compact pairs.
Thus when two compact pairs $(X,X_0)$ and $(Y,Y_0)$ have the same ANR-shape, for an arbitrary abelian group $G$, we have that
$$
\vH^*(X,X_0;\mathfrak{A})\approx\vH^*(Y,Y_0;\mathfrak{A}).
$$
We omit the coefficient group $\mathfrak{A}$ in the following and the reader can take the coefficient group as the integer group.
By considering the \v Cech cohomology groups of each compact pair $(K,K_0)\ss(X,X_0)$, we can obtain an inverse system of \v Cech cohomology groups
$$
\vH^*(C(X,X_0))=\{\vH^*(K,K_0),i^*_{(K,K_0)(K',K'_0)},c(X,X_0)\}.
$$

Recall the {\em inverse limit} (see \cite{Eil}) of an inverse system of groups $G_*=(G_\lam,p_{\lam\lam'},\Lam)$ consists of a group $G_\8$ and homomorphisms $p_\lam:G_\8\ra G_\lam$ such that
\be\label{5.1}
p_{\lam\lam'}p_{\lam'}=p_\lam,\hs\lam\leq\lam'.
\ee
Moreover, if $p'_\lam:G'\ra G_\lam$ is another collection of homomorphisms with property (\ref{5.1}), then there is a unique homomorphism $g:G'\ra G_\8$ such that
$$
p_\lam g=p'_\lam,\hs\lam\in \Lam.
$$
We denote $\disp G_\8=\limla G_*$.
Clearly, the inverse limit $G_\8$ of $G_*$ is unique up to a natural isomorphism.
According to Theorem 3.14 in Chapter VIII of \cite{Eil}, the inverse limit $\disp\limla G_*$ exists for every inverse system $G_*$ of groups.
\bd
Let $(X,X_0)$ be a pair of Hausdorff spaces.
The {\bf H-shape cohomology group} $C\vH^q(X,X_0)$ for each $q\in\N$ is defined as
$$
C\vH^q(X,X_0)=\limla \vH^q(C(X,X_0)).
$$
\ed
Since \v Cech cohomology groups satisfy (P), so do the H-shape cohomology groups via the inverse limit.
Furthermore, by the properties of inverse limit, H-shape cohomology group is H-shape invariant, i.e., for each $q\in\N$,
\be\label{5.2}
{\rm Sh}_{\rm H}(X,X_0)={\rm Sh}_{\rm H}(Y,Y_0)\Ra C\vH^q(X,X_0)\approx C\vH^q(Y,Y_0).
\ee

It is trivial that H-shape cohomology groups satisfy Eilenberg-Steenrod Axioms except the exactness property.
But in our situation, we only need the exactness property for pairs of Hausdorff spaces with H-shape of compact pairs.
This is confirmed by the following theorem.

\bt\label{t4.2}
Suppose that $(X,X_0)$ and a compact pair have the same H-shape.
Then $C\vH^*(X,X_0)$ satisfies the exactness property.
\et

\bo
Let $(K,K_0)$ is a compact pair.
Assume that $(X,X_0)$ and $(K,K_0)$ have the same H-shape.

First we prove that $C\vH^*(K,K_0)$ satisfies the exactness property.
Since $(K,K_0)$ is compact, we have another direct system $K^*=\{(K,K_0),1_K\}$ besides $C(K,K_0)$, where the index set is a singleton and $1_K$ is the identity shape map.
It is a simple result in \cite{Rub} that $C(K,K_0)$ and $K^*$ have the same homotopy type and so do $\vH^q(C(K,K_0))$ and $\vH^q(K^*)$ as inverse systems of groups for $q\in\N$.
Thus by considering the inverse limit, we have
\be\label{5.3}
C\vH^q(K,K_0)\approx\limla\vH^q(K^*)=\vH^q(K,K_0).
\ee
As known in \cite{Eil}, \v Cech cohomology groups $\vH^*(K,K_0)$ satisfy the exactness property.
And hence by the isomorphism (\ref{5.3}), $C\vH^*(K,K_0)$ also has the exactness property.

By the supposition and (\ref{5.2}), we obtain the exactness of $C\vH^*(X,X_0)$.
The proof is complete.
\eo

Now we consider a local semiflow $\Phi$ on a complete metric space $X$.
Let $K$ be a compact isolated invariant set of $\Phi$.

\bd
Let $K$ be a compact isolated invariant set with $(N,E)$ being a shape index pair of $K$.
The {\bf H-shape cohomology index} of $K$ is defined for each $q\in\N$, as
$$
C\vH^q(s(\Phi,K))=C\vH^q(N/E,[E]).
$$
\ed
If the semiflow $\Phi$ is clear, we simply write the H-shape cohomology index of $K$ as $C\vH^*(s(K))$.

As H-shape cohomology groups satisfy the property (P) and (\ref{5.2}), by Theorem \ref{t:4.1}, it is easy to see that $C\vH^q(s(K))$ is independent of the choice of shape index pairs.
Therefore, H-shape cohomology index is well defined as above.

\subsection{Morse Equations}\label{ss5.3}

Morse equation is one of the most interesting and important topics of compact invariant sets for dynamical systems (\cite{Con1,Kap,Rob,Ryb,Sanj1,Wang}).
It can reflect a lot of information of the inner structure of compact invariant set, such as the dimensions, the topological structure of the Morse sets and the connected trajectories between Morse sets.
We establish the Morse equations in the framework of H-shape cohomology index in a standard way as follows.

Let $X$ be a complete metric space and $K$ a compact isolated invariant set.
Suppose $K$ has a Morse decomposition $\sM=\{M_1,\,\cdots,\,M_n\}$ with the corresponding Morse filtration
$$
\es=A_0\ss A_1\ss\cdots\ss A_n=K.
$$
Let $(N,E)$ be a shape index pair of $K$. Let $\~\Phi$ be the quotient flow on $N/E$ defined as in Section \ref{s3}.
Then it can be shown that $\~\sM=\{\~M_0,\~M_1,\,\cdots,\,\~M_n\}$ forms a Morse decomposition of the global attractor $\cA$ of $\~\Phi$, where
$$
\~M_0=\{[E]\},\hs \~M_k=\pi(M_k)\,\,(1\leq k\leq n),
$$
and $\pi:N\cup E\ra N/E$ is the quotient map.
Then we have a corresponding Morse filtration $\{\es,\cA_0,\cdots,\cA_n\}$ of $\~\sM$.
Moreover $\cA_k$ is closed for each $k=0,\cdots,n$, since $N/E$ is a normal Hausdorff space.

By Theorem \ref{th3.6}, for $k=0$, 1, $\cdots$, $n-1$, $\cA_k$ has a $\cK_0^\8$ Lyapunov function $\zeta_k$ on the region of attraction $\W_k:\,=\W(\cA_k)$.
Pick $a_k>0$ and set
$$
\~N_k:\,=\Cap_{i=k}^{n-1}\zeta_i^{a_i},\hs k=0,\,1,\,\cdots,\,n-1.
$$
Let $N_n=N$ and $N_k=\pi^{-1}(\~N_k)$, $k=0$, 1,$\cdots$, $n-1$. Then we have a sequence of closed subsets satisfying
$$
N_0\ss N_1\ss\cdots \ss N_n.
$$

It is easy to verify that $(N_k,N_{k-1})$ and $(N_k,N_0)$ are shape index pairs of $M_k$ and $A_k$, respectively, $k=1$, $\cdots$, $n$.
By very standard argument (see e.g. \cite{Ryb,Sanj1}) one can obtain the following {\em Morse equation} associated with H-shape cohomology theory:
\be\label{e5.5}
\sum_{k=1}^{n}\sum_{q=0}^\8t^q\,{\rm rank}C\vH^q(N_k,N_{k-1})=\sum_{q=0}^\8t^q\,{\rm rank}C\vH^q(N_n,N_0)+(1+t)Q(t),
\ee
where
$$
Q(t)=\sum_{k=1}^{n}\sum_{q=1}^{\8}t^{q-1} {\rm rank}\,\de_{k}^q
$$
and $\de_{k}^q$ is the coboundary operator from $C\vH^{q-1}(N_{k-1},N_0)$ to $C\vH^q(N_k,N_{k-1})$.

Referring to the property (P), we know for $0\leq l\leq k\leq n$, $C\check H^q(N_k,N_l)$ is isomorphic to $C\vH^q(N_k/N_l,[N_l])$.
As $(N_k,N_{k-1})$ and $(N_k,N_0)$ are shape index pairs of $M_k$ and $A_k$, respectively, we have
$$
C\vH^q(N_k,N_{k-1})\approx C\vH^q(s(M_k)), \hs C\vH^q(N_k,N_0)\approx C\vH^q(s(A_k)).
$$
Hence \eqref{e5.5} can be rewritten as follows:
\be\label{e5.4}
\sum_{k=1}^{n}\sum_{q=0}^\8t^q\,{\rm rank}C\vH^q(s(M_k))=\sum_{q=0}^\8t^q\,{\rm rank}C\vH^q(s(K))+(1+t)Q(t).
\ee

For each compact isolated invariant set $M$, set
$$
p(t,s(M))=\sum_{q=0}^{\8}t^q\,{\rm rank}C\check{H}^q(s(M)).
$$
Then $p(t,\,s(M))$ is called  the {\em formal Poincar\'{e} polynomial} of $s(M)$.
Now the Morse equation \eqref{e5.4} can be restated in terms of formal Poincar\'{e} polynomials,
$$
\sum_{k=1}^{n}p(t,s(M_{k}))=p(t,s(K))+(1+t)Q(t).
$$
\br We consider the calculation of the Morse equation of the maximal compact invariant set for a general nonlinear evolutionary equation
\be\label{6.1}\frac{\di u}{\di t}+\cN u+f(u)=0,\Hs u\in X,\ee
where $X$ is a Banach space, $\cN:X\ra X^*$ is a nonlinear operator and $f:X\ra X^*$ is a nonlinear functional.

In a large number of applications (see, \cite{Hen,Ryb,Tem}), the maximal compact invariant set $K$ of \eqref{6.1} is contained in some (finite dimensional) subspace of $X$.
Indeed, Morse equation of $K$ only relies on the unstable manifold of $K$, but not the spaces containing $K$.
Thus shape index pairs' independence of neighborhoods for isolated invariant sets brings great convenience to the relevant calculations.

The topological structure of $K$ is usually very complicated, let alone the dynamical behaviors of the restricted systems on $K$.
Besides, the natural Lyapunov function of (\ref{6.1}) (whenever existing) is often defined on some proper subspace or even some lower-dimensional invariant manifold in $X$.
As a result, one can hardly give the Conley index pairs of the Morse set of $K$ via merely the Lyapunov function of (\ref{6.1}).
However, our shape index pairs can well avoid such difficulties.

The example (Section 6) in \cite{Wang} illustrates the flexible calculations of our shape index pairs and Morse theory.
Moreover, the theory in this paper can applied to the case when the phase space is not separable.
\er

\section{Applications to a Retarded Nonautonomous System}

\subsection{Basic Notations and Results}
Let $\cH$ be a compact metric space with metric $\di(\cdot,\cdot)$.
A given dynamical system $\theta$ is defined on $\cH$, i.e., a  continuous mapping $\theta:\R\X\cH\ra\cH $ satisfying the following group property:
$$
\theta_0h=h,\Hs \theta_{s+t}h=\theta_s\theta_th
$$
for all $h\in\cH$ and $s,t\in\R$. (Here we have written $\theta(t,h)$ as $\theta_th$.)

We assume $\cH$ is {\it minimal} for $\theta$, that is, the dynamical system $\theta$ has no nonempty compact invariant proper subsets of $\cH$.
\vs

Let $X$ be a real Banach space with norm $\|\cdot\|$ and $L$ be a sectorial operator on $X$ with compact resolvent.
Pick a number $a>0$ such that
\be\label{7.1}
\hbox{Re}\,z\geq a_0>0,\Hs\mb{for all }\,z\in\sig(L+aI)\mb{ and a constant }a_0.
\ee
Set $\cL=L+aI$.
Then we can define the fractional powers $\cL^\al$ for $\al\in(0,1)$; see \cite{Hen} for details.
For each $\al\geq0$, define the fractional power $X^\al$ of $X$ to be the space $D(\cL^\al)$, which is equipped {with the norm} $\|\cdot\|_\al$ defined by
$$\|x\|_\al=\|\cL^\al x\|,\hs x\in X^\al.$$
Note that the definition of $X^\al$ is independent of the choice of the number $a$.
We denote by $C_\al$ the embedding constant from $X^\al$ into $X$, i.e., $\|x\|\leq C_\al\|x\|_\al$, for $x\in X^\al$.

Suppose that $L$ has a spectral decomposition $\sig(L)=\sig^-\cup\sig^+$, where
\be\label{7.2}
\mb{Re}\,z\leq -\b<0\,\,(z\in\sig^-),\hs \mb{Re}\,z\geq \b>0\,\,(z\in\sig^+)
\ee
for some $\b>0$.
Let $X=X_1\oplus X_2$ be the corresponding direct sum decomposition of $X$ with $X_1$ and $X_2$ being invariant subspaces of $L$.
Let $\Pi_k:X\ra X_k$ ($k=1,2$)  be the projection from $X$ to $X_k$.
Denote $L_k=L|_{X_k}$ and by $\{e^{-L_kt}\}_{t\geq0}$ the semigroup generated by $-L_k$.
By the basic knowledge on sectorial operators (see Henry \cite{Hen}),  we know that there exists $M\geq1$ such that
\be\label{7.3}\ba{c}\|\cL^\al e^{-L_1t}\|\leq Me^{\b t},\hs \|e^{-L_1t}\|\leq Me^{\b t},\Hs t\leq0,\\
\|\cL^\al e^{-L_2t}\Pi_2\cL^{-\al}\|\leq Me^{-\b t},\hs\|\cL^\al e^{-L_2t}\|\leq Mt^{-\al}e^{-\b t},\Hs t>0.\ea\ee

\subsection{Main Problem and Conclusion}
Consider the following retarded  cocycle system in $X$:
\be\label{7.4}
\left\{\ba{lr}\disp\frac{\di u}{\di t}+Lu=f(\theta_th,u,u(\cdot-\tau)),&t>0,\,h\in \cH,\\[2ex]
u(t)=\vsig(t),&t\in[-\tau,0],\ea\right.
\ee
where $f:\cH\X X^\al\X X\ra X$ for some $\al\in[0,1)$ and $\vsig:[-\tau,0]\ra X$ are continuous and the nonnegative number $\tau$ is the time delay.
The space $\cH$ and the system $\theta$ are usually called the {\em base space} and the {\em driving system} of (\ref{7.4}), respectively. We make the following assumption.\\
\noindent{\em (F1)}\hs The nonlinear term $f(h,x,y)$ is globally Lipschitzian in $(x,y)$ in a uniform manner with respect to $h\in\cH$, namely, there exists $l>0$ such that
$$
\|f(h,x,y)-f(h,x',y')\|\leq l(\|x-x'\|_\al+\|y-y'\|)
$$
for all $h\in \cH$, $x,x'\in X^\al$ and $y,y'\in X$.

Under this assumption, we infer from the basic theory on retarded evolution equations in Banach spaces (see \cite{LiW,WuJ}) that,
the system (\ref{7.4}) has a unique global mild solution $u:[-\tau,\8)\ra X$ for each initial data $\vsig\in\cC:=C([-\tau,0],X)$.
Here a {\em global mild solution} $u=u(t;\vsig,h)$ for (\ref{7.4}) is a mapping $u$ from $[-\tau,\8)\X\cC\X\cH$ to $X$ such that $u(t;\vsig,h)\in X^\al=D(\cL^\al)$ for $t>0$, $\vsig\in\cC$, $h\in\cH$
and $u$ satisfies the integral problem
\be\label{7.5}
\left\{\ba{l}\disp u(t)=e^{-Lt}\vsig(0)+\int_{0}^te^{-L(t-\nu)}f(\theta_\nu h,u(\nu),u(\nu-\tau))\di\nu,\;t>0,\\[2ex]
u(t)=\vsig(t),\;-\tau\leq t\leq0.\ea\right.
\ee

Furthermore, we assume that $f$ satisfies the following condition.\\
{\em (F2)}\hs The nonlinear term $f$ is sublinear at the infinity in the subspace $X^\al\X X^\al$ uniformly on $\cH$, i.e.,
$$
\frac{\|f(h,x,y)\|}{\|x\|_\al+\|y\|_\al}\ra0\hs\mb{ as }\|x\|_\al+\|y\|_\al\ra\8,\mb{ for all }h\in\cH.
$$

Then we have the existence of bounded full solutions for (\ref{7.4}) as follows.
\bt\label{t7.1}
Suppose $\cH$ is minimal for $\theta$.
Under the assumptions  (F1) and (F2) on $f$, the system (\ref{7.4}) has a bounded full solution pertaining to each $h\in\cH$,
i.e., for each $h\in\cH$, there exists a bounded full solution $u(t;\vsig,h)$ defined for a certain $\vsig\in\cC$ and all $t\in\R$.
\et

As a generalisation of a similar result for nonautonomous systems in \cite{JuLi},
this theorem can be verified by using Banach contraction mapping principle in complete metric spaces.
In our situation, we will employ H-shape index to prove this theorem in the following subsection.

\subsection{The Proof of Theorem \ref{t7.1}}

First we consider the following equation dependent on $\lam\in[0,1]$,
\be\label{7.7}
\left\{\ba{l}\disp\frac{\di u}{\di t}+Lu=\lam f(\theta_th,u,u(t-\tau)),\Hs t>0,\,h\in \cH,\\[2ex]
u(t)=\vsig(t),\;-\tau\leq t\leq0,\ea\right.
\ee
where $\lam\in[0,1]$.
Note that when $\lam=0$, (\ref{7.7}) is a linear equation; when $\lam=1$, (\ref{7.7}) is our original equation (\ref{7.4}).
Let $\vp_\lam:[-\tau,\8)\X\cH\X\cC\ra X$ be the solution mapping of (\ref{7.7}) for each $\lam\in[0,1]$.
Then by the classical theory of retarded, nonautonomous functional differential equations (see \cite{Hen,Rob,WuJ}),
$\vp_\lam(t;h,\vsig)$ is continuous in $\lam,t,h,\vsig$ respectively.

Concerning the equation (\ref{7.7}), we have the following result.

\bl\label{l7.2}
There exists $R>0$ such that, for every bounded full solution $u:\R\ra X$ of (\ref{7.7}) with each $\lam\in[0,1]$, we have that
\be\label{7.8}
u(t)\in X^\al\hs\mb{and}\hs\|u(t)\|_\al\leq R\Hs\mb{for all }t\in\R.
\ee
\el
\bo
This conclusion is a generalisation of a relevant one for autonomous dynamical systems (Theorem 5.1 in Chapter II of \cite{Ryb}).

By the theory of functional differential equations (see \cite{Hen,WuJ}), if $u$ is a full solution of (\ref{7.7}) in $X$, $u(t)$ is in $X^\al$ for all $t\in\R$.

We prove (\ref{7.8}) by contradiction.
Suppose that for every $R>0$, there exists a full solution $u(t)$ of (\ref{7.7}) such that $\|u(t)\|_\al>R$ for some $t\in\R$ and $\lam\in[0,1]$.
Then there is a sequence $\lam_n\in[0,1]$ and a sequence of full bounded solutions $t\ra u_n(t)$ of (\ref{7.7}) with $\lam=\lam_n$ such that
\be\label{7.9}
c_n:=\sup_{t\in\R}\|u_n(t)\|_\al\ra\8,\hs\mb{ as }n\ra\8,
\ee
and $\|u_n(t_n)\|_\al>c_n-1$ for some $t_n\in\R$.

Let $v_n=c_n^{-1}u_n$ and $f_n:\cH\X X^\al\X X\ra X$ be defined as
\be\label{7.10}
f_n(h,x,y)=c_n^{-1}\lam_nf(h,c_nx,c_ny).
\ee
We claim that for each $R_0\geq0$,
\be\label{7.11}
S_n:=\sup\{\|f_n(h,x,y)\|:\|x\|_\al+\|y\|_\al\leq R_0,h\in\cH\}\ra0\hs\mb{ as }n\ra\8.
\ee

Assume that this claim holds.
Noting that $\|v_n(t)\|_\al$ is bounded by $1$ for $t\in\R$, we take $R_0\geq2$ in (\ref{7.11}).
By \cite{JuLi}, the bounded full solution $u(t)$ of (\ref{7.4}) satisfies the following integral equation,
\be\label{7.12}
\ba{ll}u(t)=&\disp\int_{-\8}^te^{-L_2(t-\nu)}\Pi_2f(\theta_\nu h,u(\nu),u(\nu-\tau))\di\nu\\[2ex]
&\disp-\int_t^{\8}e^{-L_1(t-\nu)}\Pi_1f(\theta_\nu h,u(\nu),u(\nu-\tau))\di\nu.\ea
\ee
By substitution of $u,\,f$ by $v_n,\,f_n$, respectively, (\ref{7.12}) also satisfies (\ref{7.7}).
As a result, by the inequalities (\ref{7.3}),
$$\|v_n(t_n)\|_\al\leq S_nM\(\int_{-\8}^{t_n}(-\nu)^{-\al}e^{\b\nu}\di\nu+\int_{t_n}^\8e^{-\b\nu}\di\nu\)\ra0,\hs\mb{as }n\ra\8.$$
However, by (\ref{7.9}),
$$1\geq\|v_n(t_n)\|_\al=c_n^{-1}\|u_n(t_n)\|_\al>\frac{c_n-1}{c_n}\ra1,\hs\mb{as }n\ra\8,$$
which leads to a contradiction! This asserts the lemma.
\vs

Now it remains to show the claim (\ref{7.11}).
Indeed, by the assumption (F2), for every $\ve>0$, there is $R_1>0$ such that for all $h\in\cH$,
$$
\|f(h,x,y)\|\leq\ve(\|x\|_\al+\|y\|_\al)\hs\mb{ for }\|x\|_\al+\|y\|_\al>R_1.
$$
By the assumption (F1), for $\|x\|_\al+\|y\|_\al\leq R_1$, we have
$$\|f(h,x,y)\|\leq l(\|x\|_\al+C_\al\|y\|_\al)+\|f(h,0,0)\|\leq l(1+C_\al)R_1+m,$$ where $m=\max_{h\in\cH}\|f(h,0,0)\|$.
When $\|x\|_\al+\|y\|_\al\leq R_0$, we obtain from (\ref{7.10}) that
$$
\|f_n(h,x,y)\|\leq\left\{\ba{ll}c_n^{-1}[l(1+C_\al)R_0+m],&\mb{ if }c_n(\|x\|_\al+\|y\|_\al)\leq R_1,\\
\ve R_0,&\mb{ if }c_n(\|x\|_\al+\|y\|_\al)>R_1.\ea\right.
$$
This implies the claim. The proof is thus complete.
\eo

For sake of applying H-shape index, we take the space $\cC$ into consideration as the phase space for the system generated by (\ref{7.7}) as follows (see \cite{WuJ}).

We endow a norm $\|\cdot\|_\cC$ on $\cC$,  such that $\|\vsig\|_\cC=\max_{-\tau\leq s\leq0}\|\vsig(s)\|$.
For two real numbers $t'\leq t''$, $t\in[t',t'']$ and a continuous function $u:[t'-\tau,t'']\ra X$, we denote by $u_t$ the element of $\cC$ given by $u_t(s)=u(t+s)$ for $t\in[t',t'']$ and $s\in[-\tau,0]$.
Similar to (\ref{7.5}), the mild solution of (\ref{7.7}) can be written as
\be\label{7.13}
\left\{\ba{ll}u_t(s)=&\disp e^{-L(t+s)}\vsig(0)+\lam\int_0^{t+s}e^{-L(t+s-\nu)}f(\theta_\nu h,u_\nu(0),u_\nu(-\tau))\di\nu,\\[2ex]
&t>0,\,-\tau\leq s\leq0,\\
u_0=\vsig.&\ea\right.
\ee
In this framework, we denote the solution mapping of (\ref{7.7}) by $\~\vp_\lam:\R^+\X\cH\X\cC\ra\cC$ such that for $s\in[-\tau,0]$,
$$
\~\vp_\lam(t;h,\vsig)(s)=\vp_\lam(t+s;h,\vsig)=u_t(s)
$$
with $u_t(s)$ defined as (\ref{7.13}).
The continuity of $\~\vp_\lam$ over $\R^+\X\cH\X\cC$ follows immediately from that of $\vp_\lam$ and the compactness of $[-\tau,0]$ for each $\lam\in[0,1]$.
Also $\vp_\lam$ depends on $\lam\in[0,1]$ continuously.
This allows us to define a skew product flow $\Phi_\lam$ on $\cH\X\cC$ for each $\lam\in[0,1]$ as
\be\label{7.13A}
\Phi_\lam(t)(h,\vsig)=(\theta_th,\~\vp_\lam(t;h,\vsig)),
\ee
where we endow the space $\cH\X\cC$ a metric $\~\di(\cdot,\cdot)$ such that
$$
\~\di((h_1,\vsig_1),(h_2,\vsig_2))=\di(h_1,h_2)+\|\vsig_1-\vsig_2\|_\cC.
$$
Furthermore, we define a subspace $\cC^\al$ of $\cC$ such that
$$\cC^\al=\{\vsig\in\cC:\vsig(t)\in X^\al,\mb{ for all }t\in[-\tau,0]\},$$
with the norm $\|\vsig\|_{\cC^\al}=\max_{-\tau\leq s\leq0}\|\vsig(s)\|_\al$ for every $\vsig\in\cC^\al$.

Let $\sK_\lam$ be the union of all full bounded orbits of $\Phi_\lam$ in $\cH\X\cC$ for each $\lam\in[0,1]$.
By Lemma \ref{l7.2}, there is $R>0$ such that
\be\label{7.14}
\sK_\lam\ss\cH\X\mB^\al(R)\ss\cH\X\mB(C_\al R),\Hs\mb{for all }\lam\in[0,1],
\ee
where $\mB^\al(R)$ and $\mB(R)$ denote the open balls centred at the origin ${\bf 0}$ with radius $R\geq0$ in $\cC^\al$ and $\cC$, respectively.
By the property of mild solutions, the unstable manifold of $\sK_\lam$ is surely contained in $\cH\X\cC^\al$.
Thus, we have the following conclusion.
\bl\label{l7.3}
The pair $(\Phi_\lam,\sK_\lam)$ is S-continuous on $[0,1]$.
\el
\bo
Since $f$ is globally Lipschitzian on $X^\al\X X$ uniformly in $h$, by the classical results (\cite{Hen,Ryb,WuJ}), we know that $\Phi_\lam$ does not explode in every bounded set of $\cH\X\cC$ for $\lam\in[0,1]$.
Then by the definition of S-continuity and the continuity of $\Phi_\lam$ in $\lam$, we only need to verify that the product set $\cH\X\ol\mB(C_\al R)$ is $\{\Phi_{\lam_n}\}$-admissible for every sequence $\lam_n\in[0,1]$ with $\lam_n\ra\lam$.

Let $h_n\in\cH$, $\vsig_n\in\ol\mB(C_\al R)$, $0<t_n\ra\8$ such that $\Phi_{\lam_n}([0,t_n])\vsig_n\ss\ol\mB(C_\al R)$.
We will show that $\Phi_{\lam_n}(t_n)(h_n,\vsig_n)$ has a convergent subsequence in $\cH\X\cC$.
By the definition of global mild solutions, we can assume $t_n>2\tau$ for all $n\in\N^+$ and hence we have that $\~\vp_{\lam_n}(t_n,h_n,\vsig_n)\in\cC^\al$.
Note that $\cH$ is compact.
By Arzela-Ascoli theorem, it is sufficient to show that the sequence $\vp_{\lam_n}(t_n+\cdot;h_n,\vsig_n)$ is equi-continuous on $[-\tau,0]$ in $X$ for all  $n$, i.e., for every $\ve>0$, there exists $\de>0$ such that for all $s_1,\,s_2\in[-\tau,0]$ with $|s_1-s_2|<\de$,
$$
\|\vp_{\lam_n}(t_n+s_1;h_n,\vsig_n)-\vp_{\lam_n}(t_n+s_2;h_n,\vsig_n)\|<\ve\hs\mb{for all }n\in\N^+.
$$

Indeed, denoting $u_t^n(s)=\vp_{\lam_n}(t+s;h_n,\vsig_n)$ and applying the representation (\ref{7.13}), we have for $-\tau\leq s_1\leq s_2\leq0$,
\begin{align}
&u_{t_n}^n(s_2)-u_{t_n}^n(s_1)\notag\\
=&[e^{-L(s_2-s_1)}-I]u_{t_n}^n(s_1)+\lam_n\int_{t_n+s_1}^{t_n+s_2}e^{-L(t_n+s_2-\nu)}f(\theta_\nu h_n,u_\nu^n(0),u_\nu^n(-\tau))\di\nu\notag\\
:=&J_1+J_2,\label{7.14A}\end{align}
where $I:X\ra X$ is the identity.

It is deduced from Theorem 1.4.3 in \cite{Hen} that, if $s\geq0$ and $u\in X^\al$ with arbitrary $\al\in(0,1)$, then there is $\kappa_{\al}>0$ such that
$$
\|(e^{-\cL s}-I)u\|\leq\al^{-1}\kappa_{\al}s^\al\|u\|_\al.
$$
Hence we have
$$
\ba{ll}&\|(e^{-Ls}-I)u\|\\\leq&\|(e^{-\cL s}-I)e^{as}u\|+\|(e^{as}-1)u\|\\
\leq&\al^{-1}\kappa_{\al}s^{\al}e^{as}\|u\|_\al+C_\al(e^{as}-1)\|u\|_\al\\
\leq&Cs^\al\|u\|_\al,\ea
$$
where the positive constant $C$ depends only on $\al$.
The last inequality holds for all sufficiently small $s<1$.
Thus, we have that when $s_2-s_1$ is sufficiently small,
\be\label{7.16}
\|J_1\|\leq C(s_2-s_1)^\al\|u_{t_n}^n(s_1)\|_\al.
\ee

By the assumption (F1) and the discussion in the proof of Lemma \ref{l7.2}, we concludes that for a certain $m>0$,
$$
\|f(h,x,y)\|\leq m+lC_\al R+l\|x\|_\al,\mb{ for }\|y\|\leq C_\al R.
$$
It follows from (\ref{7.1}) that $\|e^{-Lt}\|\leq e^{\b_0t}$, for all $t\geq0$ and a certain $\b_0>0$.
And hence
\be\label{7.17}
\ba{ll}\|J_2\|&\disp\leq \int_{s_1}^{s_2}\|e^{-L(s_2-\nu)}f(\theta_\nu h_n,u_{t_n}^n(\nu),u_{t_n}^n(\nu-\tau))\|\di\nu\\[2ex]
&\disp\leq\int_{s_1}^{s_2}e^{\b_0(s_2-\nu)}\(m+l C_\al R+l\|u_{t_n}^n(\nu)\|_\al\)\di\nu.\ea
\ee

By the estimates \eqref{7.3}, there are $M_0,\,\b_0>0$ such that for all $t\geq0$,
$$\|\cL^\al e^{-Lt}\|\leq M_0t^{-\al}e^{\b_0t}.$$
Hence for $s\in[-\tau,0]$,
\begin{align}&\|u_{t_n}^n(s)\|_{\al}\notag\\
\leq&\|e^{-L(s+2\tau)}u_{t_n-2\tau}^n(0)\|_\al+\int_{-2\tau}^s\|e^{-L(s-\nu)}f(\theta_\nu h_n,u^n_{t_n}(\nu),u^n_{t_n}(\nu-\tau))\|_\al\di\nu\notag\\
\leq&M_0(s+2\tau)^{-\al}e^{\b_0(s+2\tau)}\|u_{t_n-2\tau}^n(0)\|\notag\\
&+M_0\int_{-2\tau}^s(s-\nu)^{-\al}e^{\b_0(s-\nu)}\(m+lC_\al R+l\|u_{t_n}^n(\nu)\|_\al\)\di\nu\notag\\
\leq&M_0e^{2\b_0\tau}\left[\tau^{-\al}C_\al R+\int_{-2\tau}^s(s-\nu)^{-\al}\(m+lC_\al R+l\|u_{t_n}^n(\nu)\|_\al\)\di\nu\right]\notag\\
\leq&M_0e^{2\b_0\tau}\left[\tau^{-\al}C_\al R+(m+lC_\al R)\frac{(2\tau)^{1-\al}}{1-\al}+l\int_{-2\tau}^s(s-\nu)^{-\al}\|u_{t_n}^n(\nu)\|_\al\di\nu\right]\notag
\end{align}
By Gronwall's inequality (see Lemma 7.1.1 in \cite{Hen}), we easily know that
\be\label{7.17A}\|u_{t_n}^n(s)\|_{\al}\leq M',\ee
for $s\in[-\tau,0]$, where $M'=M'(\tau,\al,l,m,R)$ and is independent of the sequence $\lam_n$.

Combining (\ref{7.14A}), (\ref{7.16}), (\ref{7.17}) and (\ref{7.17A}), we obtain the desired equi-continuity for the sequence $u_t^n$.
This makes $\Phi_{\lam_n}(t_n)(h_n,\vsig_n)$ possess a convergent subsequence.
Hence $\cH\X\ol\mB(C_\al R)$ is strongly $\{\Phi_{\lam_n}\}$-admissible, which indicates that $(\Phi_\lam,\sK_\lam)$ is S-continuous on $[0,1]$.
The proof is complete.
\eo

Since $L$ has a compact resolvent, according to (\ref{7.1}) and (\ref{7.2}), the set $\sig^-$ contains only finitely many eigenvalues.
Let $r$ be the sum of all multiplicities of the eigenvalues in $\sig^-$.
It is known that $r>0$.

For the calculation of H-shape index $s(\Phi_1,\sK_1)$, we also need to notice the following simple result,
$$
s(\Phi_1\X\Phi_2,K_1\X K_2)=s(\Phi_1,K_1)\wedge s(\Phi_2,K_2),
$$
for two disjoint semiflows $\Phi_1$ and $\Phi_2$ with their compact isolated invariant sets $K_1$ and $K_2$, respectively, where
\be\label{4.1}(\Phi_1\X\Phi_2)(t)(x_1,x_2):=(\Phi_1(t)x_1,\Phi_2(t)x_2).\ee

Now we calculate the H-shape index as follows.
\bl\label{l7.4} The H-shape index of $\sK_1$ for $\Phi_1$ is
$$s(\Phi_1,\sK_1)=\cH^*\wedge\Sig^r\ne\ol0,$$
where $\cH^*$ denotes the H-shape of the pointed space $(\cH\cup\{*\},*)$ with $*\notin\cH$.
\el
\bo We split the proof  into two steps.
\vs

\noindent{\it Step 1.} We first consider the linear equation, that is, (\ref{7.7}) when $\lam=0$,
$$
\frac{\di u}{\di t}+Lu=0,\;t>0,\Hs u(s)=\vsig(s),\;-\tau\leq s\leq0.
$$

In consideration of the phase space $\cC$ and recalling the consequences in \cite{WuJ}, we obtain that the function $w(t):=u_t$  satisfies the linear system,
\be\label{7.18}
\frac{\di w}{\di t}+L_ww=0,\;t>0,\Hs w(0)=\vsig,
\ee
where $L_w$ is a sectorial operator on $\cC$ corresponding to $L$ on $X$.
Furthermore, the operator $L_w$ has the same eigenvalues as $L$ with the same multiplicities, respectively.

We denote by $\phi$ the semiflow generated by the linear equation (\ref{7.18}).
Observe that the origin ${\bf 0}$ of $\cC$ is the maximal compact invariant set of $\phi$ in $\cC$, and ${\bf 0}$ is a hyperbolic point for (\ref{7.18}).
Therefore, the unstable subspace of ${\bf 0}$ in $\cC$ is $r$-dimensional.
By this fact, it is easy to obtain that
$$s(\phi,\{{\bf 0}\})=\Sig^r.$$

\noindent{\it Step 2.} We consider the continuous family of semiflows $\Phi_\lam$, $\lam\in[0,1]$, defined in (\ref{7.13A}), and compute the H-shape index of $\sK_1$ for $\Phi_1$ via the continuation property.

When $\lam=0$, we see that $\Phi_0=\theta\X\phi$, which is defined as (\ref{4.1}).
By the results of Step 1, we have that for every $R>0$,
$$
s(\Phi_0,\sK_0)=s(\Phi_0,\cH\X\{{\bf 0}\})=s(\theta,\cH)\wedge s(\phi,\{{\bf 0}\})=\cH^*\wedge\Sig^r.
$$

By the S-continuity of $(\Phi_\lam,\sK_\lam)$ at $\lam\in[0,1]$ stated in Lemma \ref{l7.3}, we infer from Theorem \ref{th4.9} that
$$
s(\Phi_1,\sK_1)=s(\Phi_0,\sK_0)=\cH^*\wedge\Sig^r.
$$
Since $\cH\cup\{*\}$ is not connected, it is clear that $\cH^*\wedge\Sig^r\ne\ol0$.
The calculation is finished.
\eo
\br Here we adopt H-shape index instead of Conley index, due to the flexible choices of the index pairs for H-shape index.
We only need to consider the unstable subspace (manifold), no matter how complicated the original phase space is.
\er
In the following, we use the framework of $\Phi=\Phi_1:\cH\X\cC\ra\cH\X\cC$ and denote simply
$$\Phi(t)(h,\vsig)=(\theta_th,\~\vp(t;h,\vsig))\hs\mb{and}\hs\sK=\sK_1.$$
Based on Lemma \ref{l7.2}, \ref{l7.3} and \ref{l7.4}, we now prove the main result Theorem \ref{t7.1}.
\Vs

\noindent{\it Proof of Theorem \ref{t7.1}.}
By Lemma \ref{l7.4}, Example \ref{x4.4} and Lemma \ref{l7.2}, we know that $\sK\ne\es$ and $\sK\ss\cH\X\ol\mB^\al(R)$ for some $R>0$.

We claim that for each $h\in\cH$, there is $\vsig\in\cC$ such that $(h,\vsig)\in\sK$.

Suppose that this claim is true.
Noting that $\sK$ is an invariant set of $\Phi$.
For each $h\in\cH$, there is a full solution $\~\gam_h$ of $\Phi$ contained in $\sK$ such that
$$\~\gam_h(t)=(\theta_th,u^h_t),\hs\mb{for all }t\in\R,$$
with $u^h_t$satisfying $u^h_t=\~\vp(t-t';\theta_{t'}h,u^h_{t'})$ for all $t,t'\in\R$ and $t\geq t'$.
Converting the phase space from $\cC$ back to $X$, we have a full solution $u^h$ of (\ref{7.4}) such that
$$u^h(t)=u^h_t(0).$$
Then $u^h$ is a full solution of (\ref{7.4}) pertaining to $h$ with $\|u^h(t)\|\leq\|u^h_t\|_\cC\leq C_\al R$ for all $t\in\R$.
This leads to the final conclusion of Theorem \ref{t7.1}.

Now it remains to prove the claim.
Since the driving system $\theta$ is independent of the phase space, the projection $P:\cH\X\cC\ra\cH$ and the systems $\Phi$, $\theta$ satisfy the following commutativity,
$$\theta_t\circ P=P\Phi(t).$$
Because $\sK$ is invariant for $\Phi$, we have
$$\theta_t(P\sK)=P\Phi(t)\sK=P\sK,\hs\mb{for all }t\geq0,$$
which implies $P\sK$ is invariant for $\theta$.
Moreover, by the compactness of $\sK$, we know that $P\sK$ is compact in $\cH$.
Therefore, $P\sK$ is a compact invariant set for $\theta$ in $\cH$.

By the minimality of $\cH$ for $\theta$, the compact invariant sets in $\cH$ are only $\es$ and $\cH$ itself.
Whereas, $\sK$ is nonempty, and so is $P\sK$.
As a result, we must have that $P\sK=\cH$, which implies the claim.
\qed

\bibliographystyle{plain}

\begin{thebibliography}{99}

\small\setlength{\baselineskip}{10pt} \vskip 10pt
\bibitem{Bor}Borsuk K., {\it Theory of shape}. Monografie Matematyczne, Tom 59. PWN-Polish Scientific, Warsaw, 1975.
\bibitem{Con1}Conley C., {\it Isolated invariant sets and the Morse index}. Regional Conference Series in Mathematics, vol. 38, American Mathematical Society, Providence, RI, 1978.
\bibitem{Eil}Eilenberg S., Steenrod N., {\it Foundation of algebraic topology}. Princeton University Press, Princeton, N.J., 1952.
\bibitem{Gir1}Giraldo A., Mor\'{o}n M.A., Ruiz del Portal F.R., Sanjurjo J.M.R., Shape of global attractors in topological spaces. {\it Nonlinear Analysis}, 60(5): 837 - 847, 2005.
\bibitem{Gir2}Giraldo A., Jim\'{e}nez R., Mor\'{o}n M.A., et al. Pointed shape and global attractors for metrisable spaces. {\it Topology and its Appl.}, 158(2): 167 - 176, 2011.
\bibitem{Hat}Hatcher A., {\it Algebraic Topology}. Cambridge University Press, 2002.
\bibitem{Hen}Henry D., {\it Geometric theory of semilinear parabolic equations}, Lect. Notes in Math. 840, Springer Verlag, Berlin New York, 1981.
\bibitem{JuLi}Ju X., Li D., Global synchronising behavior of evolution equations with exponentially growing nonautonomous forcing. {\it Comm. Pure Appl. Analysis}, 17(5): 1921 - 1944, 2018.
\bibitem{JuQi}Ju X., Qi A., Wang J., Strong Morse-Lyapunov functions for Morse decompositions of attractors of random dynamical systems. {\it Stochastics and Dynamics}, 18(01): 1850012, 2018.
\bibitem{Kap}Kapitanski L., Rodnianski I., Shape and Morse theory of attractors. {\it Comm. Pure Appl. Math.}, LIII: 0218 - 0242, 2000.
\bibitem{Li1}Li D., Morse theory of attractors via Lyapunov functions. {\it arXiv preprint}, arXiv: 1003.1576, 2009.
\bibitem{Li2}Li D., Shi G., Song X., A linking theory for dynamical systems with applications to PDEs. {\it arXiv preprint}, arXiv:1312.1868v3, 2015.
\bibitem{Li3}Li D., Wang J., Xiong Y., Attractors of local semiflows on topological spaces. {\it J. Korean Math. Soc.}, 54(3): 773 - 791, 2017.
\bibitem{LiW}Li D., Wei J., Wang J., On the dynamics of abstract retarded evolution equations. {\it Abstract and Applied Analysis}, 2013, 2013.
\bibitem{Mar2}Marde\v{s}i\'{c} S., Segal J., {\it Shape theory. The inverse system approach}. North-Holland Mathematical Library, 26. North-Holland, Amsterdam - New York, 1982.
\bibitem{Marz}Marzocchi A., Zandonella Necca S., Attractors for dynamical systems in topological spaces. {\it Discrete Contin. Dyn. Syst.}, 8(3): 585¨C597,  2002.
\bibitem{Mroz}Mrozek M., Shape index and other indices of Conley type for local maps on locally compact Hausdorff spaces. {\it Fundam. Math.}, 145: 15 - 37 1994.
\bibitem{Rob}Robbin J.W., Salamon D., Dynamical systems, shape theory and the Conley index. {\it Ergodic Theory Dynam. Syst.}, 8: 375 - 393, 1988.
\bibitem{Rub}Rubin L., Sanders J., Compactly generated shape. {\it Gen. Topology Appl.}, 4: 73 - 83, 1974.
\bibitem{Ryb}Rybakowski K.P., {\it The homotopy index and partial differential equations}. Springer-Verlag, etc. 1980.
\bibitem{Sanc2}S\'{a}nchez-Gabites J.J., An approach to the Conly shape index without index pairs. {\it Rev. Mat. Complut.}, 24: 95 - 114, 2011.
\bibitem{Sand1}Sanders T.J., Shape groups for Hausdorff spaces. {\it Glasnik Mat. ser. iii}, 297 - 304, 1973.
\bibitem{Sand2}Sanders T.J., Shape groups and products. {\it Pacific Journal of Mathematics}, 48(2): 485-496, 1973.
\bibitem{Sanj1}Sanjurjo Jos\'{e} M.R., Morse equation and unstable manifolds of isolated invariant sets. {\it Nonlinearity}, 16: 1435 - 1448, 2003.
\bibitem{Sanj2}Sanjurjo Jos\'{e} M.R., Shape and Conley index of attractors and isolated invariant sets. {\it Differential Equations, Chaos and Variational Problems}. Birkh\"auser Basel, 2008.
\bibitem{Tem}Temam R., {\it Infinite-Dimensional Dynamical Systems in Mechanics and Physics} (Second Edition). Springer-Verlag, New York, 1997.
\bibitem{Wang}Wang J., Li D., Duan J., On the shape Conley index theory of semiflows on complete metric spaces. {\it Discrete and Continuous Dynamical Systems - Series A}, 36(3): 1629 - 1647, 2015.
\bibitem{WuJ}Wu J., {\it Theory and applications of partial functional differential equations}, Springer-Verlag, New York, Inc., 1996.
\end{thebibliography}

\end{document}